\documentclass[10pt,a4paper]{amsart}
\usepackage[T1]{fontenc}
\usepackage[british]{babel}
\usepackage{amssymb,mathrsfs}
\IfFileExists{newtxtext.sty}{\usepackage{newtxtext,newtxmath}}{\usepackage{mathptmx}}
\usepackage{needspace}
\usepackage{microtype}
\usepackage[margin=1in,heightrounded]{geometry}
\usepackage[hidelinks]{hyperref}
\hypersetup{pdftitle={Singularities in leaf closures on complex surfaces},pdfauthor={}}

\newtheorem{theorem}{Theorem}
\newtheorem{proposition}[theorem]{Proposition}
\newtheorem{lemma}[theorem]{Lemma}
\newtheorem{corollary}[theorem]{Corollary}
\theoremstyle{remark}
\newtheorem{remark}[theorem]{Remark}

\newcommand{\C}{\mathbb C}
\newcommand{\PP}{\mathbb P}
\newcommand{\J}{\mathscr J}
\DeclareMathOperator{\Sing}{Sing}
\DeclareMathOperator{\med}{med}
\DeclareMathOperator{\Fol}{Fol}

\title[Singularities in leaf closures]
{Singularities in leaf closures on complex surfaces}

\author{Mauricio Corr\^ea}
\address{(M. Corr\^ea) Department of Mathematics,
	University of Bari Aldo Moro,
	Via Edoardo Orabona 4,
	70125 Bari, Italy}
\email{mauricio.barros@uniba.it}

\author{Fernando Reis}
\address{(F. Reis) Polytechnic Institute,
	Federal University of Rio de Janeiro,
	Avenida Alu\'izio da Silva Gomes, 50,
	Bairro da Gl\'oria,
	27930-560 Maca\'e, RJ, Brazil}
\email{fernandoreis@ipoli.macae.ufrj.br}


\begin{document}
	
	\begin{abstract}
		We prove a geometric and dynamical criterion excluding minimal
		sets contained in the regular locus of a holomorphic foliation
		on a compact complex surface. In particular, we prove that every
		Jouanolou foliation of degree at least two has no exceptional
		minimal set. We also give an explicit lower bound, growing
		linearly with the degree, for the number of singularities in
		every regular leaf closure of a Jouanolou foliation. This bound
		persists under small perturbations, as does the existence of a
		nonempty open saturated set whose leaves accumulate at all
		singularities. Under finite pullback unramified over the singular
		set, the bound is multiplied by the degree of the normalised
		component over the base leaf. Together with a general persistence
		theorem, this yields nonempty analytic open sets of foliations
		on every smooth projective surface with arbitrarily large
		uniform lower bounds on the number of singularities in their
		regular leaf closures. On double-plane K3 surfaces, we obtain
		foliations with ample canonical bundle and explicitly describe
		and count their singularities.
	\end{abstract}
	
	\maketitle
	
	\section{Introduction}
	
	We study conditions ensuring that the closure of every regular
	leaf of a holomorphic foliation on a compact complex surface
	contains a prescribed minimum number of singular points. The
	problem includes the question of exceptional minimal sets on
	\(\PP^2\): \textit{does every regular leaf accumulate at a singularity?}
	Camacho, Lins Neto and Sad proved that there is at most one
	exceptional minimal set and that it meets every algebraic
	curve \cite{CLS}. Bonatti, Langevin and Moussu proved that
	such a set contains a leaf with contracting hyperbolic
	holonomy \cite{BLM}.
	
	Our approach begins with domains in which real curves tangent
	to the foliation converge to a singular point. Entry into
	these domains is proved from finite jets along regular leaves,
	or from holonomy maps which decrease a nonnegative function
	on compact transverse sets. A curve with Stein complement
	then connects these local constructions to every leaf closure.
	When several curves are used, intersection estimates among
	the resulting finite sets of singularities give a stronger
	bound. The precise criterion is
	Theorem~\ref{thm:propagation-count}; it allows several
	attracting domains at the same singularity and counts each
	singular point only once.
	
	When a curve with Stein complement is available, bounds obtained
	from hyperbolic singularities persist. Here hyperbolic means that the two eigenvalues
	are nonzero and their ratio is not real. For a holomorphic line
	bundle \(N\), let \(\Fol(X,N)\) denote the space of foliations
	with normal bundle \(N\), represented by sections of
	\(\Omega_X^1\otimes N\) with isolated zeros, up to scalar.
	We write \(L_q\) for the regular leaf through \(q\) and use
	the analytic topology on \(\Fol(X,N)\). The next theorem is proved in
	Section~\ref{sec:persistence} by small attracting neighbourhoods,
	finite leafwise paths and compactness.
	
	\Needspace{5\baselineskip}
	\begin{theorem}\label{thm:persistence}
		Let \(\mathcal F_0\in\Fol(X,N)\), and let \(P\) be a finite
		set of hyperbolic singularities of \(\mathcal F_0\).
		Let \(D\subset X\setminus\Sing(\mathcal F_0)\) be a compact
		complex curve with Stein complement. Suppose that, for some
		integer \(r\geq1\),
		\[
		\#(\overline{L_q}\cap P)\geq r\qquad(q\in D).
		\]
		There is a neighbourhood \(\mathcal V\) of \(\mathcal F_0\)
		in \(\Fol(X,N)\) such that every regular leaf \(L\) of every
		\(\mathcal F\in\mathcal V\) satisfies
		\[
		\#\bigl(\overline L\cap P(\mathcal F)\bigr)\geq r,
		\qquad P(\mathcal F)=\{p(\mathcal F):p\in P\}.
		\]
		Furthermore, if a regular point \(q_0\) satisfies
		\(P\subset\overline{L_{q_0}}\) for \(\mathcal F_0\), there
		are an open neighbourhood \(B\) of \(q_0\) and a possibly
		smaller \(\mathcal V\) such that, for every
		\(\mathcal F\in\mathcal V\) and \(q\in B\), the closure
		of the \(\mathcal F\)-leaf through \(q\) contains
		\(P(\mathcal F)\). Its saturation is a nonempty open set
		with the same property.
	\end{theorem}
	
	On a projective surface, a general ample divisor provides the
	curve required by this theorem. Finite maps to the plane give
	examples with lower bounds as large as desired. We obtain the
	following consequence in Section~\ref{sec:finite}.
	
	\Needspace{5\baselineskip}
	\begin{corollary}\label{cor:any-surface}
		Let \(Y\) be a smooth projective complex surface and
		\(r\geq1\). There are a line bundle \(N_r\) and a nonempty
		analytic open subset \(\mathcal W_r\subset\Fol(Y,N_r)\)
		such that every regular leaf \(L\) of every
		\(\mathcal F\in\mathcal W_r\) satisfies
		\[
		\#\bigl(\overline L\cap\Sing(\mathcal F)\bigr)\geq r.
		\]
		The open set contains a foliation without invariant algebraic
		curves. No member has a nonempty compact saturated subset
		in its regular locus.
	\end{corollary}
	
	The normal bundle is fixed on each open set in this corollary
	and may vary with \(r\). We construct a pullback foliation
	without invariant algebraic curves in that open set; the
	dynamical conclusion holds for every member of the open set.
	The elementary diagonal foliations of
	Section~\ref{sec:linear-example} illustrate the local criteria
	and contracting holonomy before the higher-degree application.
	
	The main explicit application concerns the Jouanolou foliation
	\(\J_d\) on \(\PP^2\), induced by
	\[
	J_d=z_2^d\frac{\partial}{\partial z_1}
	+z_3^d\frac{\partial}{\partial z_2}
	+z_1^d\frac{\partial}{\partial z_3},\qquad d\geq2.
	\]
	Its degree is the number of tangencies with a general line.
	Jouanolou proved that it has no invariant algebraic curve
	\cite{Jouanolou}. Its regular leaves are therefore Zariski
	dense, although their closures in the usual topology can be
	proper subsets of the plane.
	
	\Needspace{5\baselineskip}
	\begin{theorem}\label{thm:main}
		For \(d\geq2\), put
		\(k_d=\lfloor(d+1)/2\rfloor\) and
		\(b_d=\max\{k_d,3k_d-3\}\). Every regular leaf \(L\)
		of \(\J_d\) satisfies
		\[
		\#\bigl(\overline L\cap\Sing(\J_d)\bigr)\geq b_d
		=\begin{cases}
			1,&d=2,\\
			3\left\lfloor\dfrac{d+1}{2}\right\rfloor-3,&d\geq3.
		\end{cases}
		\]
		In particular, \(\J_d\) has no exceptional minimal set.
	\end{theorem}
	
	Camacho and de Figueiredo proved the absence of exceptional
	minimal sets for \(2\leq d\leq5\), using interval arithmetic
	and validated integration of a real flow \cite{CF}. They also
	proved that the leaf through \([0:0:1]\) accumulates at all
	singularities in every degree. The diagonal symmetries fix
	that regular point and preserve its leaf. For an arbitrary
	leaf, a symmetry can change the leaf, so one must construct
	and transport the paths which give its singular limits.
	Alvarez and Deroin established structural stability in degree
	two \cite{AD}; Gehlawat has announced structural stability
	in every degree \cite{Gehlawat}. The broader class of
	foliations of Jouanolou type is discussed in
	\cite[Section~4]{Deroin}. Our proof of the bound and its
	persistence uses the local and compact arguments given here.
	
	The proof of Theorem~\ref{thm:main} has three parts. First,
	we construct an attracting domain \(U_d\) where the arguments
	of the three homogeneous coordinates have spread
	\(\Delta<\pi/(2d)\). Multiplying \(J_d\) by a phase
	determined by their median makes the extreme arguments move
	towards one another. The radial equations give positive lower
	bounds for the normalised moduli; these bounds make angular
	contraction exponential. A Lyapunov function then forces the
	ratios of the moduli to tend to one, giving convergence to
	\([1:1:1]\). The estimate \(|\beta_i|\leq d\Delta\)
	is responsible for the width of the sector. The exponential
	estimate refers to the flow time of the normalised real field
	(Corollary~\ref{cor:sector-rate}).
	
	Second, the first nonzero derivatives of the coordinates along
	a leaf starting on \(\{z_2=0\}\) show that it enters every
	admissible diagonal translate of \(U_d\). The admissible
	rotations form an arc containing at least \(k_d\) consecutive
	singularity indices. The other coordinate lines give blocks
	with two different steps; any two blocks from different lines
	have at most one common index. Third, the maximum principle
	forces every leaf closure to meet all three lines. Saturation
	puts the corresponding singularities in the original closure,
	and inclusion--exclusion gives \(b_d\).
	
	The same construction gives a nonempty open saturated set
	whose leaves accumulate at every singularity. This property
	and the lower bound persist under small perturbations by
	Theorem~\ref{thm:persistence}; see
	Corollary~\ref{cor:robustness}. The sector itself need not
	be preserved under perturbation, since the paths can be
	continued into smaller local attracting domains.
	
	For a finite pullback, we prove that a leaf closure maps onto
	the closure of a base leaf and that the singularity count is
	multiplied by the degree of its normalised component over that
	leaf (Theorem~\ref{thm:finite-pullback}). We also compute the
	normal and canonical bundles and the total Milnor number.
	For a double plane branched over a general smooth sextic,
	the surface is K3 and the pulled-back foliation has ample
	canonical bundle. Its singularities consist of
	\(2(d^2+d+1)\) hyperbolic points and \(6(d+5)\)
	nondegenerate points with a Morse first integral
	(Theorem~\ref{thm:k3}). This gives a direct comparison between
	the singularities used in the dynamical bound and the global
	Chern number.
	
	\section{Leaf closures on compact complex surfaces}\label{sec:criterion}
	
	Throughout this section, \(X\) is a compact complex surface
	and \(\mathcal F\) is a holomorphic foliation by curves with
	finite singular set \(S\). Closures are taken in \(X\), and
	leaves carry their intrinsic topology. For a regular point \(q\),
	write \(L_q\) for its leaf. A symmetry is a holomorphic
	automorphism preserving \(\mathcal F\). A subset of
	\(X\setminus S\) is saturated if it contains the entire leaf
	through each of its points.
	
	\Needspace{5\baselineskip}
	\begin{lemma}\label{lem:closure-saturation}
		If \(L\) is a regular leaf and \(K=\overline L\), then
		\(K\setminus S\) is saturated in the regular locus.
	\end{lemma}
	
	\begin{proof}
		Let \(q\in K\setminus S\), and choose a flow box
		\(D\times T\) centred at \(q=(0,0)\), with plaques
		\(D\times\{\xi\}\). There are points
		\((u_n,\xi_n)\in L\) tending to \((0,0)\). Each plaque
		\(D\times\{\xi_n\}\) lies in \(L\); hence, for every
		\(u\in D\), the points \((u,\xi_n)\) tend to \((u,0)\).
		The plaque through \(q\) is therefore contained in \(K\).
		For any regular leaf \(M\) meeting \(K\), this proves that
		\(M\cap K\) is open in \(M\). It is also closed in \(M\).
		Connectedness gives \(M\subset K\).
	\end{proof}
	
	The next lemma applies even when the compact set already
	contains singularities. Logarithmic poles place the maximum
	at a regular point, where saturation can be used.
	
	\Needspace{5\baselineskip}
	\begin{lemma}\label{lem:saturation}
		Let \(D\subset X\) be a compact complex curve such that
		\(X\setminus D\) admits a smooth strictly plurisubharmonic
		function. If \(K\subset X\) is compact, has a regular point,
		and \(K\setminus S\) is saturated, then \(K\cap D\ne\varnothing\).
		In particular, the closure of every regular leaf meets \(D\).
	\end{lemma}
	
	\begin{proof}
		Suppose \(K\cap D=\varnothing\), and put \(V=X\setminus D\).
		Let \(\varphi\) be a smooth strictly plurisubharmonic function
		on \(V\). For each \(s\in S\cap V\), choose disjoint
		relatively compact coordinate balls \(B_s\subset V\), with
		coordinates \(w_s\) centred at \(s\). Let \(\chi_s\) be a
		smooth cutoff function supported in \(B_s\), equal to one on
		a smaller ball. On \(V\setminus S\), set
		\[
		h_s=\chi_s\log\|w_s\|^2,
		\]
		extending it by zero outside \(B_s\).
		
		Near \(s\), the function \(h_s\) is plurisubharmonic. The
		possible negative part of its Levi form is supported in a
		compact annulus away from \(s\). There are finitely many
		such annuli, and the Levi form of \(\varphi\) has a positive
		lower bound there relative to any fixed Hermitian metric.
		Consequently, for a sufficiently large constant \(A>0\),
		\[
		\rho=A\varphi+\sum_{s\in S\cap V}h_s
		\]
		is strictly plurisubharmonic on \(V\setminus S\). If
		\(S\cap V\) is empty, take \(\rho=\varphi\). At each
		\(s\in S\cap V\), extend \(\rho\) by \(-\infty\).
		Its logarithmic term tends to \(-\infty\), while the other
		terms are smooth near \(s\), so this extension is upper
		semicontinuous.
		
		The function \(\rho\) is bounded above on the compact set
		\(K\), and is finite at a regular point of \(K\). It therefore
		attains a finite maximum at some \(q\in K\setminus S\).
		The regular leaf through \(q\) lies in \(K\) by saturation.
		On a holomorphic plaque parametrised by an immersion
		\(\gamma\), with \(\gamma(0)=q\), strict plurisubharmonicity
		gives \(\Delta(\rho\circ\gamma)(0)>0\). This contradicts
		the maximum at \(q\).
	\end{proof}
	
	If \(V\) is Stein, a smooth strictly plurisubharmonic exhaustion
	is available; see \cite[Chapter~2]{Forstneric}. If \(X\) is
	projective and \(D\) is an effective ample divisor, choose
	\(m\) such that \(\mathcal O_X(mD)\) is very ample. In the
	associated projective embedding, the section with divisor
	\(mD\) defines a hyperplane. Its complement on \(X\) is a
	closed affine subvariety, so \(X\setminus D\) is Stein.
	
	For the complement of a line in \(\PP^2\), there is a direct
	formula. In affine coordinates \(w=(w_1,w_2)\), with
	\(T=S\cap\C^2\), take
	\begin{equation*}
		\rho(w)=|w_1|^2+|w_2|^2+
		\sum_{s\in T}\log\|w-s\|^2.
	\end{equation*}
	The Levi form of a logarithmic summand, evaluated at
	\(v\in\C^2\), is
	\[
	\frac{\|v\|^2\|w-s\|^2-|\langle v,w-s\rangle|^2}
	{\|w-s\|^4}\geq0
	\]
	by Cauchy--Schwarz. The squared norm makes the sum strictly
	plurisubharmonic. Its restriction to a plaque satisfies
	\[
	\Delta(\rho\circ\gamma)(0)
	\geq4\bigl(|\gamma_1'(0)|^2+|\gamma_2'(0)|^2\bigr)>0,
	\]
	which is the contradiction used in Lemma~\ref{lem:saturation}.
	
	A compact saturated subset of the regular locus contains a
	minimal such subset. On a general surface this minimal set may
	be a compact leaf; the term \emph{exceptional minimal set}
	excludes compact leaves and the whole ambient surface. Our
	conclusion excludes every nonempty compact saturated subset
	of the regular locus, and therefore excludes these possibilities
	together.
	
	\Needspace{5\baselineskip}
	\begin{lemma}\label{lem:zorn}
		The following are equivalent:
		\begin{enumerate}
			\item the closure of every regular leaf meets \(S\);
			\item no nonempty compact saturated subset is contained in
			\(X\setminus S\);
			\item no minimal compact saturated subset is contained in
			\(X\setminus S\).
		\end{enumerate}
	\end{lemma}
	
	\begin{proof}
		If a compact saturated set \(E\subset X\setminus S\)
		contains a point \(q\), it contains \(\overline{L_q}\),
		which is disjoint from \(S\). Thus (1) implies (2), and
		(2) implies (3). Conversely, if a regular leaf has closure
		\(K\) disjoint from \(S\), then \(K\) is compact and
		saturated by Lemma~\ref{lem:closure-saturation}. Order its
		nonempty compact saturated subsets by reverse inclusion.
		The intersection of a chain is nonempty by compactness and
		is again compact and saturated. Zorn's lemma supplies a
		minimal compact saturated subset of \(K\), contradicting (3).
	\end{proof}
	
	On \(\PP^2\), a compact regular leaf cannot occur. Indeed,
	such a leaf \(C\) would be a smooth algebraic curve. Flow
	boxes along \(C\) have transverse transition maps depending
	only on the transverse coordinate. Their derivatives along
	\(C\) are locally constant, so the normal bundle of \(C\)
	is flat and has degree zero. On the other hand, its degree
	is \(C^2=(\deg C)^2>0\). Thus, on \(\PP^2\), the minimal
	sets in Lemma~\ref{lem:zorn} are exactly the exceptional
	minimal sets considered in \cite{CLS,CF}.
	
	We shall use two facts about curves tangent to the foliation.
	The first follows from flow boxes; the second is the uniqueness
	statement needed when a real curve is defined by a holomorphic
	generator multiplied by a continuous complex function of time.
	
	\Needspace{5\baselineskip}
	\begin{lemma}\label{lem:tangent-curve}
		Let \(I\subset\mathbb R\) be an interval and
		\(\gamma:I\to X\) a \(C^1\) curve.
		\begin{enumerate}
			\item If \(\gamma(I)\subset X\setminus S\) and
			\(\dot\gamma(t)\) is tangent to \(\mathcal F\) at every
			\(t\), then \(\gamma(I)\) lies in a single regular leaf.
			\item In a coordinate chart, suppose that
			\(\dot\gamma(t)=\alpha(t)V(\gamma(t))\), where \(V\)
			is a holomorphic generator of \(\mathcal F\) and
			\(\alpha:I\to\C\) is continuous. If \(\gamma(t_0)\)
			is regular for some \(t_0\in I\), then \(\gamma(t)\)
			is regular for every \(t\in I\).
		\end{enumerate}
	\end{lemma}
	
	\begin{proof}
		In a flow box with transverse coordinate \(\xi\), tangency
		gives \(d\xi(\dot\gamma)=0\). The curve therefore lies
		locally in a plaque. Since \(I\) is connected, these plaques
		belong to a single leaf, proving (1).
		
		For (2), regard the chart as an open subset of \(\mathbb R^4\).
		The map \((t,z)\mapsto\alpha(t)V(z)\) is continuous in time
		and locally Lipschitz in the state variable, uniformly on
		compact time intervals. Uniqueness holds in both time
		directions. If \(\gamma(t_*)=s\) is a zero of \(V\), the
		constant curve at \(s\) solves the same equation. The set
		\(\{t\in I:\gamma(t)=s\}\) is closed by continuity and
		open by local uniqueness. It is therefore all of \(I\),
		contrary to regularity at \(t_0\).
	\end{proof}
	
	Convergence to a singular point can be checked by a Lyapunov function after the orbit enters a
	compact set. This form is useful when radial estimates provide
	compactness before convergence is known.
	
	\Needspace{5\baselineskip}
	\begin{proposition}\label{prop:lyapunov-criterion}
		Let \(U\subset X\) contain a unique singularity \(p\), and
		let \(W\) be a locally Lipschitz real vector field on \(U\),
		vanishing at \(p\) and tangent to \(\mathcal F\) on the
		regular locus. Suppose its positive flow is defined for all
		time in \(U\). For each positive orbit \(\gamma_q\) starting
		at a regular point, assume that there exist \(T\geq0\), a
		compact set \(K_q\subset U\), and a nonnegative locally
		Lipschitz function \(H_q\) near \(K_q\) such that
		\[
		\gamma_q(t)\in K_q\quad(t\geq T),\qquad
		H_q^{-1}(0)\cap K_q=\{p\},\qquad
		\frac{d}{dt}H_q(\gamma_q(t))\leq-c_qH_q(\gamma_q(t))
		\]
		for almost every \(t\geq T\), with \(c_q>0\).
		Then \(\gamma_q\) stays in \(L_q\) at every finite time
		and converges to \(p\).
		
		If, in addition, \(H_q(x)\geq a_q\operatorname{dist}(x,p)^\nu\)
		on \(K_q\), where \(a_q,\nu>0\), the convergence is
		exponential, with rate \(c_q/\nu\).
	\end{proposition}
	
	\begin{proof}
		Since \(W\) is locally Lipschitz and \(W(p)=0\), uniqueness
		prevents a regular orbit from reaching \(p\) at finite time.
		There are no other singularities in \(U\). The orbit remains
		regular and belongs to \(L_q\) by
		Lemma~\ref{lem:tangent-curve}(1).
		
		The function \(H_q\circ\gamma_q\) is locally absolutely
		continuous. Integrating its differential inequality gives
		\[
		0\leq H_q(\gamma_q(t))
		\leq H_q(\gamma_q(T))e^{-c_q(t-T)},\qquad t\geq T.
		\]
		On the compact part of \(K_q\) outside any neighbourhood
		of \(p\), the minimum of \(H_q\) is positive. Hence the
		orbit eventually lies in every such neighbourhood and
		converges to \(p\). Under the final hypothesis,
		\[
		\operatorname{dist}(\gamma_q(t),p)
		\leq\left(\frac{H_q(\gamma_q(T))}{a_q}\right)^{1/\nu}
		e^{-c_q(t-T)/\nu}.
		\]
		Enlarging the constant includes the interval \([0,T]\).
	\end{proof}
	
	For the counting argument, let \(P\subset S\) be finite.
	For each \(p\in P\), allow a finite family of open
	neighbourhoods \(U_{p,\beta}\), each containing no other
	singularity. On every such domain assume there is a locally
	Lipschitz real field \(W_{p,\beta}\), tangent to the foliation,
	vanishing at \(p\), with complete positive flow converging
	to \(p\). For a regular point \(q\), define
	\[
	A(q)=\{p\in P:L_q\cap U_{p,\beta}\ne\varnothing
	\text{ for some }\beta\}.
	\]
	Every point of \(A(q)\) lies in \(\overline{L_q}\).
	
	\Needspace{5\baselineskip}
	\begin{proposition}\label{prop:counting}
		Let \(D_1,\ldots,D_m\subset X\setminus S\) be compact
		complex curves with Stein complements. Suppose that, for each
		\(q\in D_i\), a set \(B_i(q)\subset A(q)\) of cardinality
		at least \(r_i\geq1\) is specified. If
		\[
		\#\bigl(B_i(q)\cap B_j(q')\bigr)\leq b_{ij}
		\quad(q\in D_i,\ q'\in D_j,\ i<j),
		\]
		then every regular leaf \(L\) satisfies
		\begin{equation}\label{eq:general-count}
			\#(\overline L\cap S)\geq
			\max\left\{\max_i r_i,\ \sum_i r_i-\sum_{i<j}b_{ij}\right\}.
		\end{equation}
		For a single curve, the hypothesis \(\#A(q)\geq r\) on
		that curve gives \(\#(\overline L\cap S)\geq r\).
	\end{proposition}
	
	\begin{proof}
		Put \(K=\overline L\). Lemmas~\ref{lem:closure-saturation}
		and \ref{lem:saturation} give \(q_i\in K\cap D_i\) and
		\(L_{q_i}\subset K\). Following the positive flow in any
		domain met by \(L_{q_i}\) shows that
		\(B_i(q_i)\subset K\cap S\). Their union contains each
		of the sets, hence has cardinality at least \(\max_i r_i\).
		The first two terms of inclusion--exclusion give
		\[
		\#\bigcup_iB_i(q_i)
		\geq\sum_i\#B_i(q_i)
		-\sum_{i<j}\#\bigl(B_i(q_i)\cap B_j(q_j)\bigr)
		\geq\sum_i r_i-\sum_{i<j}b_{ij}.
		\]
		Indeed, a point belonging to \(a\geq1\) sets contributes
		\(a-\binom a2\leq1\) to the middle expression. The
		single-curve assertion is the case \(m=1\).
	\end{proof}
	
	\section{Entry into attracting domains}
	
	We retain the attracting domains \(U_{p,\beta}\) and fields
	\(W_{p,\beta}\) introduced after
	Proposition~\ref{prop:lyapunov-criterion}. When a single domain
	is considered, write it as \(U_p\) and its field as \(W_p\).
	We give two sufficient conditions for a regular leaf to meet
	such a domain.
	
	\subsection{Finite jets at a boundary point}
	
	The first condition concerns the leading coefficients of local
	functions along a leaf. The functions may be local equations of
	components of a divisor; they need not be coordinates. Write
	\(\mathbb S^1=\{z\in\C:|z|=1\}\).
	
	\Needspace{5\baselineskip}
	\begin{lemma}\label{lem:jet-entry}
		Let \(q\in X\setminus S\), and let \(V\) be a nonvanishing
		holomorphic generator of \(\mathcal F\) near \(q\), with local
		complex flow \(\Phi_s\). Let \(f_1,\ldots,f_\ell\) be holomorphic
		functions on a neighbourhood \(B\) of \(q\). Suppose that, for an
		open set \(\Theta\subset(\mathbb S^1)^\ell\),
		\[
		\Omega=\left\{x\in B:f_\nu(x)\ne0\ \text{for all }\nu,
		\ \left(\frac{f_\nu(x)}{|f_\nu(x)|}\right)_{\nu=1}^\ell
		\in\Theta\right\}\subset U_p.
		\]
		For each \(\nu\), suppose there is a finite integer \(m_\nu\geq0\)
		such that \(V^j f_\nu(q)=0\) for \(j<m_\nu\) and
		\(V^{m_\nu}f_\nu(q)\ne0\). Here \(V^0f_\nu=f_\nu\).
		If some \(c\in\C^*\) satisfies
		\[
		a_\nu=\frac{c^{m_\nu}}{m_\nu!}V^{m_\nu}f_\nu(q),
		\qquad
		\left(\frac{a_\nu}{|a_\nu|}\right)_{\nu=1}^\ell\in\Theta,
		\]
		then \(\Phi_{ct}(q)\in U_p\) for all sufficiently small real
		\(t>0\). In particular, \(p\in\overline{L_q}\).
		Moreover, a neighbourhood of \(q\) consists of points whose
		leaves meet \(U_p\).
	\end{lemma}
	
	\begin{proof}
		Taylor expansion in complex time gives
		\(f_\nu(\Phi_{ct}(q))=a_\nu t^{m_\nu}+O(t^{m_\nu+1})\).
		For real \(t>0\), the factor \(t^{m_\nu}\) is positive.
		Thus all the functions are nonzero for small \(t>0\), and their
		normalised values tend to \(a_\nu/|a_\nu|\). Since \(\Theta\)
		is open, the trajectory enters \(\Omega\). It is a regular
		curve in \(L_q\), and its subsequent \(W_p\)-orbit converges
		to \(p\) in the same leaf. To prove the last assertion, fix one
		such time \(t_0>0\). The local map \(x\mapsto\Phi_{ct_0}(x)\)
		is defined near \(q\), and its value at \(q\) belongs to the
		open set \(U_p\). Continuity gives the assertion after shrinking
		the neighbourhood.
	\end{proof}
	
	The condition uses only finitely many derivatives at the point,
	together with a specified angular region contained in an attracting
	domain. Its existence is independent of the choice of generator.
	Indeed, replacing \(V\) by \(bV\), with \(b(q)\ne0\), multiplies
	the first nonzero derivative of order \(m_\nu\) by
	\(b(q)^{m_\nu}\); replacing \(c\) by \(c/b(q)\) compensates for
	this change. If a divisor equation is multiplied by a
	holomorphic unit, its phase can vary with the point; the
	angular region is therefore not transformed by a fixed
	rotation on the whole neighbourhood. The strict condition
	on the leading coefficients permits a smaller angular region
	and a smaller neighbourhood, on which continuity of the unit
	preserves the required inclusion in \(U_p\). A function vanishing
	identically on the leaf has no such finite order; this includes
	the obstruction caused by an invariant boundary component.
	
	Points already in \(U_p\setminus S\) satisfy the condition by
	taking \(B\subset U_p\), \(f_1=1\), and
	\(\Theta=\mathbb S^1\). At a boundary point the content is the
	strict angular condition on the first nonzero derivatives.
	
	\subsection{A decreasing function for holonomy}
	
	Let \(T\) be a finite disjoint union of local transverse discs in
	\(X\setminus S\). A holonomy map between open subsets of these
	discs is obtained by following a fixed chain of flow boxes along
	a regular leafwise path. In particular, \(x\) and \(h(x)\)
	belong to the same global leaf. The domains below are actual
	domains of definition, rather than formal germs.
	
	\Needspace{5\baselineskip}
	\begin{proposition}\label{prop:holonomy-entry}
		Let \(K\subset T\) be compact, let \(E\subset K\) be relatively
		open, and suppose every point of \(E\) satisfies the finite jet
		condition of Lemma~\ref{lem:jet-entry} for \(U_p\). Let
		\(H:K\to[0,\infty)\) be continuous, with
		\(H^{-1}(0)\subset E\). Suppose that there are finitely many
		holonomy maps \(h_1,\ldots,h_s\) and a constant \(0<\lambda<1\)
		such that, for every \(x\in K\setminus E\), at least one of these
		maps is defined at \(x\) and satisfies
		\begin{equation}\label{eq:holonomy-decrease}
			h_a(x)\in K,\qquad H(h_a(x))\leq\lambda H(x).
		\end{equation}
		Then the leaf through every point of \(K\) meets \(U_p\).
		If \(K\setminus E\ne\varnothing\), put
		\(M=\max_K H\) and \(\delta=\min_{K\setminus E}H>0\).
		Every point reaches \(E\) after at most \(n\) successive
		choices in \eqref{eq:holonomy-decrease}, for any integer \(n\)
		satisfying \(\lambda^nM<\delta\).
	\end{proposition}
	
	\begin{proof}
		Start at \(x_0\in K\). If it does not belong to \(E\), choose
		a map as in \eqref{eq:holonomy-decrease} and set
		\(x_1=h_a(x_0)\). Continue while the points remain outside
		\(E\). At every defined step, \(x_j\) lies in \(K\cap L_{x_0}\)
		and \(H(x_j)\leq\lambda^jM\). The set \(K\setminus E\) is
		compact and contains no zero of \(H\), so its minimum \(\delta\)
		is positive. Remaining outside \(E\) through step \(n\) would
		give \(\delta\leq H(x_n)\leq\lambda^nM<\delta\), a contradiction.
		At the first point in \(E\), Lemma~\ref{lem:jet-entry} gives
		entry into \(U_p\). All the holonomy steps and the final local
		trajectory lie in \(L_{x_0}\).
	\end{proof}
	
	This gives a condition on defined maps, their images and a scalar
	inequality. It does not require prior knowledge of a leaf closure.
	Both the image condition and the decrease are essential: a
	contracting germ which cannot be continued along the proposed
	itinerary does not give \eqref{eq:holonomy-decrease}. Rich
	contracting holonomy dynamics also occur in the study of analytic
	continuation of holonomy maps \cite{CDFG}; the statement above
	only uses the specified finite domains and inequalities.
	
	\Needspace{5\baselineskip}
	\begin{corollary}\label{cor:contracting-holonomy}
		Let \(K=\{|z|\leq R\}\) be a closed disc in a local transversal,
		and let \(h\) be a holonomy map defined on a neighbourhood of
		\(K\), with \(h(0)=0\) and \(\sup_K|h'|\leq\kappa<1\).
		If \(0\) satisfies Lemma~\ref{lem:jet-entry} for \(U_p\), then
		the leaf through every point of \(K\) meets \(U_p\).
	\end{corollary}
	
	\begin{proof}
		Integration along the line segment from \(0\) to \(z\) gives
		\(|h(z)|\leq\kappa|z|\), so \(h(K)\subset K\). The last
		assertion of Lemma~\ref{lem:jet-entry} supplies a neighbourhood
		\(N\) of \(0\) and a fixed regular leafwise path, varying in
		flow boxes, from each point of \(N\) into \(U_p\).
		Hence \(h^n(z)\to0\) enters \(N\) after finitely many steps.
		Following this path completes the entry into \(U_p\).
	\end{proof}
	
	Here the final paths can themselves be included among the finite
	holonomy maps in Proposition~\ref{prop:holonomy-entry}. More
	generally, a holonomy map sending a compact transversal into a
	transversal contained in \(U_p\) fits that proposition: take the
	disjoint union of the source and image compact sets, let \(H\)
	be one on the source and zero on the image, and take \(E\) to
	be the image component. This also permits finite chains of
	explicit local flow maps.
	
	\subsection{Counting singularities}
	
	\Needspace{5\baselineskip}
	\begin{theorem}\label{thm:propagation-count}
		Use the attracting domains \(U_{p,\beta}\) just described.
		Let \(D_1,\ldots,D_m\subset X\setminus S\) be compact
		complex curves with Stein complements. For each \(i\),
		let \(\{O_{i\alpha}\}_\alpha\) be a finite relatively
		open cover of \(D_i\), and specify nonempty sets
		\(P_{i\alpha}\subset P\). For each \(p\in P_{i\alpha}\),
		choose one domain \(U_{p,\beta}\) for which one of the
		following holds:
		\begin{enumerate}
			\item every \(q\in O_{i\alpha}\) satisfies the finite jet
			condition of Lemma~\ref{lem:jet-entry} for that domain;
			\item \(O_{i\alpha}\) lies in a flow box whose projection
			maps it into one compact transverse set satisfying
			Proposition~\ref{prop:holonomy-entry} for that domain.
		\end{enumerate}
		Then every regular leaf \(L\) satisfies
		\begin{equation}\label{eq:target-union-bound}
			\#(\overline L\cap S)\geq
			R:=\min_{\alpha_1,\ldots,\alpha_m}
			\#\bigcup_iP_{i\alpha_i}.
		\end{equation}
		In particular, if \(\#P_{i\alpha}\geq r_i\) and
		\(\#(P_{i\alpha}\cap P_{j\alpha'})\leq b_{ij}\) for
		\(i<j\) and every \(\alpha,\alpha'\), the lower bound
		in \eqref{eq:general-count} holds. There is no nonempty
		compact saturated subset of the regular locus.
	\end{theorem}
	
	\begin{proof}[Proof of Theorem~\ref{thm:propagation-count}]
		Let \(L\) be a regular leaf and put \(K_L=\overline L\).
		Lemmas~\ref{lem:closure-saturation} and \ref{lem:saturation}
		give points \(q_i\in K_L\cap D_i\). Choose a member
		\(O_{i\alpha_i}\) of the specified cover containing \(q_i\).
		For each \(p\in P_{i\alpha_i}\), the first alternative gives
		\(L_{q_i}\cap U_{p,\beta}\ne\varnothing\) by
		Lemma~\ref{lem:jet-entry}. In the second alternative, the
		flow-box projection of \(q_i\) lies in the same leaf, and
		Proposition~\ref{prop:holonomy-entry} gives the same conclusion.
		The positive \(W_{p,\beta}\)-orbit then converges to \(p\) in
		\(L_{q_i}\). Since \(L_{q_i}\subset K_L\), we have
		\(P_{i\alpha_i}\subset K_L\cap S\).
		
		Thus \(\#(K_L\cap S)\geq\#\bigcup_iP_{i\alpha_i}\geq R\),
		which proves \eqref{eq:target-union-bound}.
		The numerical estimate follows from the first two terms of
		inclusion--exclusion, as in Proposition~\ref{prop:counting}.
		No nonempty compact saturated subset can lie in the regular
		locus when \(R\geq1\), by Lemma~\ref{lem:zorn}.
	\end{proof}
	
	The open sets in this theorem are subsets of the divisors, but
	the divisors themselves need not be transversals. At a regular
	tangency, the projection to a transverse disc is still defined
	on a flow box. The finite jet alternative likewise allows higher
	orders of contact. For a compact divisor, pointwise strict entry
	conditions give a finite cover: fix a small positive entry time
	at each point, use continuity in a neighbourhood, and take a
	finite subcover. These fixed-time paths give finite holonomy
	maps to transverse discs inside the respective \(U_p\).
	
	\Needspace{5\baselineskip}
	\begin{corollary}\label{cor:symmetry-count}
		Let a finite group \(G\) preserve \(\mathcal F\), and let
		\(U\) carry an attracting field \(W\) at \(p_0\). Use all domains
		\(g(U)\), with fields \(g_*W\), labelled by their limit
		\(g(p_0)\). They satisfy the attracting-domain assumptions
		of Theorem~\ref{thm:propagation-count}. Distinct singular
		points are indexed by \(G/G_{p_0}\), where
		\(G_{p_0}=\{g:g(p_0)=p_0\}\). Several orbits are treated
		by the disjoint union of their coset spaces.
	\end{corollary}
	
	\begin{proof}
		Pushforward preserves tangency, completeness of the positive
		flow and convergence. The map \(gG_{p_0}\mapsto g(p_0)\)
		is a bijection onto the orbit. All domains with the same
		limit are allowed in the theorem, and that point is counted
		only once.
	\end{proof}
	
	If the stabiliser does not preserve \(U\), its elements may
	give different domains at the same singularity. When a
	calculation is made at \(g^{-1}q\), push the constructed
	path forward by \(g\). Its image lies in \(L_q\) and
	converges to \(g(p_0)\), regardless of whether \(g\)
	preserves \(L_q\).
	
	\section{Persistence under perturbations}\label{sec:persistence}
	
	Let \(N\) be a holomorphic line bundle on a compact complex
	surface \(X\). We write \(\Fol(X,N)\) for the open subset of
	\(\PP H^0(X,\Omega_X^1\otimes N)\) consisting of sections
	with isolated zeros. Such a section defines a foliation with
	normal bundle \(N\); on a surface its integrability is automatic.
	All neighbourhoods in this section refer to the usual analytic
	topology on this parameter space. Local generators can be chosen
	holomorphically in the parameter after fixing local trivialisations
	and a representative in a projective coordinate chart.
	
	A singularity is called hyperbolic if its linear part has two
	nonzero eigenvalues whose ratio is not real. The following local
	fact provides attracting domains which remain available under
	perturbation.
	
	\Needspace{5\baselineskip}
	\begin{lemma}\label{lem:persistent-attractor}
		Let \(p\) be a hyperbolic singularity of \(\mathcal F_0\in\Fol(X,N)\).
		For \(\mathcal F\) near \(\mathcal F_0\), the singularity has a
		unique continuation \(p(\mathcal F)\) in a fixed coordinate
		neighbourhood. There are a constant \(c\in\C^*\), local
		holomorphic generators \(V_{\mathcal F}\), and a positive
		Hermitian form \(Q\) such that, for some \(\varepsilon,a>0\),
		the domains
		\[
		U_{\mathcal F}=\{z:(z-p(\mathcal F))^*Q(z-p(\mathcal F))
		<\varepsilon\}
		\]
		carry complete positive real-time flows of \(cV_{\mathcal F}\),
		and every orbit in \(U_{\mathcal F}\) converges to
		\(p(\mathcal F)\). The constants and the form can be chosen
		uniformly, with
		\[
		\frac{d}{dt}\bigl(w(t)^*Qw(t)\bigr)
		\leq-a\,w(t)^*Qw(t),\qquad w=z-p(\mathcal F).
		\]
	\end{lemma}
	
	\begin{proof}
		The nonreal ratio of the two eigenvalues allows a rotation
		placing both in the open left half-plane. Choose \(c\) with
		this property, and let \(A=D(cV_{\mathcal F_0})(p)\).
		The implicit function theorem gives \(p(\mathcal F)\), since
		the singularity is nondegenerate. Define
		\[
		Q=\int_0^\infty e^{A^*t}e^{At}\,dt.
		\]
		The integral converges, \(Q\) is positive definite, and
		integration of its derivative gives \(A^*Q+QA=-I\).
		In centred coordinates the nearby fields have the form
		\(cV_{\mathcal F}(p(\mathcal F)+w)=A_{\mathcal F}w+R_{\mathcal F}(w)\),
		where \(A_{\mathcal F}\to A\) and
		\(\|R_{\mathcal F}(w)\|\leq C\|w\|^2\) uniformly on a
		sufficiently small fixed ball. After shrinking the parameter
		neighbourhood,
		\(A_{\mathcal F}^*Q+QA_{\mathcal F}\leq-\tfrac12I\).
		Hence, on a smaller ball,
		\[
		\frac{d}{dt}(w^*Qw)
		\leq-\tfrac12\|w\|^2+2\|Q\|C\|w\|^3
		\leq-\tfrac14\|w\|^2
		\leq-\frac{w^*Qw}{4\|Q\|}.
		\]
		Choose \(\varepsilon\) so that the closed ellipsoids lie in
		this ball. The strict inequality on their boundaries gives
		positive invariance. Solutions remain in a compact subset of
		the coordinate neighbourhood and therefore exist for all
		positive time. The differential inequality proves convergence.
		Uniqueness prevents a regular orbit from reaching the zero
		in finite time, so the orbit remains in its original leaf.
	\end{proof}
	
	We also use the persistence of a finite path in a regular leaf.
	If such a path joins \(q\) to an interior point of
	\(U_{\mathcal F_0}\), then every point sufficiently close to
	\(q\) has a path in its leaf for \(\mathcal F\) ending in
	\(U_{\mathcal F}\), for all sufficiently close \(\mathcal F\).
	Indeed, cover the path by finitely many regular flow boxes and
	choose local transverse discs in their overlaps. Holomorphic
	flow maps for local generators give the successive holonomy
	maps. Continuous dependence on the initial point and the
	coefficients preserves their domains after shrinking the discs.
	The final point stays in \(U_{\mathcal F}\) because the original
	endpoint lies in its interior. This argument preserves fixed
	finite paths, without requiring their defining jets to retain
	the same orders of vanishing.
	
	\begin{proof}[Proof of Theorem~\ref{thm:persistence}]
		Choose disjoint coordinate neighbourhoods of the points of
		\(P\) and the attracting domains of
		Lemma~\ref{lem:persistent-attractor}. For each \(q\in D\),
		choose \(r\) distinct points \(p\in\overline{L_q}\cap P\).
		The leaf meets every chosen attracting domain, since that
		domain is a neighbourhood of its singular point. Connectedness
		of the regular leaf gives a finite path from \(q\) to an
		interior point of each domain. By the preceding flow-box
		argument, these paths persist on a neighbourhood \(B_q\)
		of \(q\), for all parameters in a neighbourhood
		\(\mathcal V_q\) of \(\mathcal F_0\).
		
		Choose finitely many \(B_q\) covering \(D\), and intersect
		their parameter neighbourhoods. Shrink the intersection so
		that \(D\) remains disjoint from the singular set; this is
		possible because the defining section is nonzero on the compact
		curve. Every point of \(D\) then has a leaf whose closure
		contains at least \(r\) distinct points of \(P(\mathcal F)\).
		For an arbitrary regular leaf \(L\),
		Lemmas~\ref{lem:closure-saturation} and \ref{lem:saturation}
		give \(q\in\overline L\cap D\) and
		\(\overline{L_q}\subset\overline L\). This proves the bound.
		
		For the final assertion, choose the finite paths from \(q_0\)
		to the domains of all points of \(P\) and take a common
		neighbourhood \(B\) on which they persist. It remains regular
		after shrinking \(\mathcal V\). The saturation of an open
		subset of the regular locus is open: locally, saturation is
		obtained by taking the union of plaques meeting the set, and
		successive flow boxes preserve openness. It is nonempty
		because it contains \(B\), and every leaf in it has the
		asserted singular limits.
	\end{proof}
	
	The theorem also applies when \(\mathcal F_0\) has other,
	possibly degenerate singularities. Only the points of \(P\)
	are continued and counted. On a projective surface a general
	member of a sufficiently ample linear system provides the
	curve \(D\). Thus any uniform lower bound obtained from
	hyperbolic singularities persists on that fixed surface and
	with that fixed normal bundle.
	
	\section{An elementary family with contracting holonomy}\label{sec:linear-example}
	
	The following degree-one example verifies the two entry
	procedures independently of the Jouanolou foliations of degree
	at least two. Its singularity count is elementary and is
	included to make the sufficient conditions explicit.
	
	\Needspace{5\baselineskip}
	\begin{proposition}
		Let \(a_0,a_1,a_2\in\C\) be noncollinear as points of
		\(\mathbb R^2\). The projective vector field induced by
		\(Z=\sum_{i=0}^2a_i z_i\partial_{z_i}\) defines a foliation
		\(\mathcal D_a\) whose singularities are the three coordinate
		points \(e_0,e_1,e_2\). For every regular point \(q=[z]\),
		\[
		\overline{L_q}\cap\Sing(\mathcal D_a)
		=\{e_i:z_i\ne0\}.
		\]
		In particular, every regular leaf accumulates at at least two
		singularities, and the uniform bound two is sharp. The general
		propagation theorem applies using any line avoiding the three
		coordinate points. For \((a_0,a_1,a_2)=(1,\omega,\omega^2)\),
		where \(\omega=e^{2\pi\mathrm i/3}\), the cyclic coordinate
		permutation is a symmetry transitive on the singular set.
	\end{proposition}
	
	\begin{proof}
		The numbers \(a_i\) are distinct, so \(Z(z)\in\C z\) holds
		exactly at the coordinate points. Each \(a_i\) is a vertex
		of their nondegenerate triangle. A real supporting linear
		functional, strictly maximised at that vertex, has the form
		\(z\mapsto\operatorname{Re}(c_i z)\), with \(c_i\ne0\).
		Consequently there is \(\gamma_i>0\) such that
		\(\operatorname{Re}(c_i(a_j-a_i))\leq-\gamma_i\) for \(j\ne i\).
		In the chart \(z_i\ne0\), the real-time flow of \(c_i Z\)
		is explicitly
		\[
		w_j(t)=w_j(0)\exp\bigl(c_i(a_j-a_i)t\bigr),\qquad
		|w_j(t)|\leq e^{-\gamma_i t}|w_j(0)|\quad(j\ne i).
		\]
		It remains in the original regular leaf at every finite time
		and converges to \(e_i\). Thus \(z_i\ne0\) implies
		\(e_i\in\overline{L_q}\). If \(z_i=0\), the invariant
		coordinate line \(\{z_i=0\}\) contains the whole leaf and
		its closure, and does not contain \(e_i\). This proves the
		equality. A regular point has at least two nonzero coordinates,
		and the regular part of a coordinate line is a leaf, proving
		sharpness.
		
		To verify the hypotheses with small attracting domains, choose
		\(U_i=\{\max_{j\ne i}|w_j|<\varepsilon\}\) in each chart.
		They carry the complete contracting flows just described. If
		\(z_i(q)\ne0\), any sufficiently large finite time sends
		\(q\) into \(U_i\). The same holds on a neighbourhood of
		\(q\). Restricting the finite-time flow to transverse discs
		gives a holonomy map whose image lies in \(U_i\); the
		two-component construction following
		Corollary~\ref{cor:contracting-holonomy} verifies
		Proposition~\ref{prop:holonomy-entry}. On a line avoiding the
		three singularities choose, near every point, two indices
		whose coordinates stay nonzero. Compactness gives the finite
		cover in Theorem~\ref{thm:propagation-count}, with one divisor
		and two singular points on every member. Its complement is affine.
		
		Finally, for the displayed triple of cube roots, cyclic
		permutation changes the diagonal field by a nonzero scalar
		after pushforward. It therefore preserves the projective
		foliation and cyclically permutes its singularities.
	\end{proof}
	
	There are also genuine contracting return maps in this example.
	On the coordinate line joining \(e_i\) to \(e_j\), choose a
	regular point with \(w_j=t_0\ne0\), and use \(u=w_k\) as
	transverse coordinate, where \(k\) is the remaining index.
	The loop \(w_j=t_0e^{2\pi\mathrm is}\), \(0\leq s\leq1\),
	has holonomy
	\[
	h(u)=\exp\left(2\pi\mathrm i\frac{a_k-a_i}{a_j-a_i}\right)u.
	\]
	Noncollinearity makes the quotient nonreal. Reversing the
	loop if necessary gives \(|h'(0)|<1\), with the same
	contraction factor on a sufficiently small closed transverse
	disc. The fixed point has explicit finite-time paths into
	both \(U_i\) and \(U_j\). Iteration of \(h\), followed by
	these paths, proves entry for every point of that disc.
	This verifies the domain and contraction conditions directly.
	
	For these linear foliations the contracting directions and
	return maps are explicit. For the higher-degree family below,
	the attracting domain will instead be constructed by choosing
	a phase from the three coordinate arguments, and finite jets
	along the coordinate lines will prove entry into it.
	
	\section{Singularities and symmetries of the Jouanolou foliation}
	
	Let \(\pi\colon\C^3\setminus\{0\}\to\PP^2\) be the
	projectivisation map, with \(\pi(z)=[z_1:z_2:z_3]\).
	Since \(\ker D\pi_z=\C z\) and \(J_d\) is homogeneous,
	\[
	D\pi_{cz}J_d(cz)=c^{d-1}D\pi_zJ_d(z)
	\qquad(c\in\C^*).
	\]
	Thus \([z]\) is singular precisely when \(J_d(z)\in\C z\).
	There are no singular points on \(\{z_3=0\}\): the third and
	first components of \(J_d(z)=\lambda z\) would successively give
	\(z_1=z_2=0\).
	
	In the chart \(z_3\ne0\), set
	\(w=(w_1,w_2)=(z_1/z_3,z_2/z_3)\). Differentiating
	the affine coordinates along \(J_d\) gives
	\[
	\begin{aligned}
		D\pi_zJ_d(z)
		&=\left(
		\frac{z_2^d z_3-z_1^{d+1}}{z_3^2},\,
		\frac{z_3^{d+1}-z_2z_1^d}{z_3^2}
		\right)\\
		&=z_3^{d-1}(w_2^d-w_1^{d+1},\,1-w_1^dw_2).
	\end{aligned}
	\]
	Thus the induced foliation is generated by
	\begin{equation}\label{eq:affine-generator}
		X_d(w_1,w_2)=(w_2^d-w_1^{d+1},\,1-w_1^dw_2).
	\end{equation}
	The affine generator has polynomial degree \(d+1\), although
	\(\J_d\) has degree \(d\). At a zero of \(X_d\), the second component gives
	\(w_1^dw_2=1\), so neither coordinate vanishes. Substituting
	\(w_2=w_1^{-d}\) into the first component yields
	\[
	w_2^d=w_1^{d+1},\qquad
	w_1^{d^2+d+1}=1.
	\]
	Thus the singular set, with \(N=d^2+d+1\), is
	\begin{equation}\label{eq:singular-set}
		\Sing(\J_d)
		=\bigl\{[\zeta:\zeta^{-d}:1]:\zeta^N=1\bigr\}.
	\end{equation}
	These are \(N\) distinct points, one of which is
	\(p_*=[1:1:1]\). Throughout the paper, fix
	\[
	\xi=e^{2\pi\mathrm i/N},\qquad
	p_j=[\xi^j:\xi^{-dj}:1],\qquad j\in\mathbb Z/N\mathbb Z.
	\]
	Thus \(p_0=p_*\). With this choice, increasing the index by one
	advances the argument of the first affine coordinate by \(2\pi/N\).
	An open arc on the unit circle therefore selects consecutive
	indices modulo \(N\). This convention will be used when the
	diagonal symmetries are counted. We call the complement of
	\(\Sing(\J_d)\) the regular locus; its leaves are the regular
	leaves of \(\J_d\).
	
	We use the notions of saturation and minimal set introduced in
	Section~\ref{sec:criterion}. In particular, on \(\PP^2\),
	excluding compact saturated subsets of the regular locus is
	equivalent to excluding exceptional minimal sets.
	
	\begin{remark}\label{rem:eigenvalues}
		The Jacobian matrix of the affine field \eqref{eq:affine-generator}
		at \(w=(1,1)\), corresponding to \(p_*\), is
		\[
		DX_d(1,1)=
		\begin{pmatrix}-(d+1)&d\\-d&-1\end{pmatrix},\qquad
		\lambda_\pm=\frac{-(d+2)\pm\mathrm i\sqrt{3}\,d}{2}.
		\]
		Both eigenvalues have negative real part, and
		\[
		\frac{\lambda_+}{\lambda_-}
		=\frac{-(d+2)+\mathrm i\sqrt3\,d}
		{-(d+2)-\mathrm i\sqrt3\,d},\qquad
		\operatorname{Im}\!\left(\frac{\lambda_+}{\lambda_-}\right)
		=-\frac{2\sqrt3\,d(d+2)}{(d+2)^2+3d^2}\ne0.
		\]
		The quotient has modulus one and is not real. In particular,
		\(\lambda_+\) and \(\lambda_-\) are linearly independent over
		\(\mathbb R\). A resonance
		\(\lambda_\pm=n_+\lambda_++n_-\lambda_-\), with
		\(n_+,n_-\in\mathbb Z_{\geq0}\), would therefore force
		\((n_+,n_-)=(1,0)\) or \((0,1)\), respectively. Hence no
		resonance of total order at least two is possible. The negative
		real parts also place the eigenvalues in the Poincar\'e domain.
		Thus the singularity is nonresonant and hyperbolic in the
		terminology of holomorphic foliations on surfaces. The nonreal
		eigenvalue condition is considered in \cite{CLS}. We verify below
		that the diagonal symmetries preserve the ratio at every
		singularity.
	\end{remark}
	
	We use the cyclic and diagonal symmetries of the
	Jouanolou foliation, as in \cite[Section~2]{CF}; see also
	\cite{ADnotes} for their description in terms of the residues
	\(1,d,d^2\) modulo \(d^2+d+1\). The cyclic permutation
	\[
	A(z_1,z_2,z_3)=(z_2,z_3,z_1)
	\]
	satisfies \(J_d(Az)=AJ_d(z)\). For each \(\zeta\) with
	\(\zeta^N=1\), where \(N=d^2+d+1\), the linear map
	\[
	D_\zeta(z_1,z_2,z_3)=(\zeta z_1,\zeta^{-d}z_2,z_3)
	\]
	satisfies
	\[
	J_d(D_\zeta z)
	=(\zeta^{-d^2}z_2^d,z_3^d,\zeta^dz_1^d)
	=\zeta^d D_\zeta J_d(z),
	\]
	where \(\zeta^{-d^2}=\zeta^{d+1}\) follows from \(\zeta^N=1\).
	Both \(A\) and \(D_\zeta\) induce linear projective
	automorphisms of \(\PP^2\) that send the tangent line of
	\(\J_d\) at each regular point to the tangent line at its
	image, and preserve \(\Sing(\J_d)\) because
	\(J_d(z)\in\C z\) if and only if the same holds for the image.
	They consequently send regular leaves to regular leaves and
	commute with taking their closures. In particular, the number
	of singularities in a leaf closure is preserved under either
	symmetry.
	
	In the affine coordinates of \eqref{eq:affine-generator},
	write \(\widetilde D_\zeta(w_1,w_2)=(\zeta w_1,\zeta^{-d}w_2)\).
	The equivariance identity becomes
	\[
	X_d(\widetilde D_\zeta w)=\zeta^d\widetilde D_\zeta X_d(w).
	\]
	Differentiation at a singular point \(w\) gives
	\[
	DX_d(\widetilde D_\zeta w)
	=\zeta^d\widetilde D_\zeta DX_d(w)\widetilde D_\zeta^{-1}.
	\]
	Thus both eigenvalues are multiplied by the same nonzero
	scalar; their ratio is unchanged. This justifies the
	hyperbolicity assertion for all singularities in
	Remark~\ref{rem:eigenvalues}.
	
	By \eqref{eq:singular-set},
	\(D_\zeta(p_*)=[\zeta:\zeta^{-d}:1]\), so the group
	\(\{D_\zeta:\zeta^N=1\}\) acts simply transitively on
	\(\Sing(\J_d)\); this will be used in the counting arguments
	below.
	
	Let \(G_d\) be the group generated by \(A\) and
	\(D_\xi\). In \(\operatorname{PGL}_3(\C)\),
	\[
	A^3=1,\qquad D_\xi^N=1,\qquad
	A D_\zeta A^{-1}=D_{\zeta^{d^2}}.
	\]
	The last identity follows by cyclically permuting the diagonal
	entries and dividing them by \(\zeta\). Thus \(G_d\) is
	finite; its elements are \(A^iD_\zeta\), with
	\(i\in\{0,1,2\}\) and \(\zeta^N=1\).
	Its action on the singularity indices is
	\begin{equation}\label{eq:cyclic-indices}
		\begin{aligned}
			A(p_j)&=[\xi^{-dj}:1:\xi^j]
			=[\xi^{-(d+1)j}:\xi^{-j}:1]=p_{d^2j},\\
			A^2(p_j)&=p_{dj}.
		\end{aligned}
	\end{equation}
	The projective normalisation divides all three coordinates by
	\(\xi^j\); the index identities use
	\(d^2\equiv-(d+1)\) and \(d^3\equiv1\pmod N\).
	
	\section{A sector criterion for the Jouanolou foliation}
	
	The convergence criterion rests on a choice of complex phase
	for the homogeneous vector field. The phase is chosen so that
	the largest and smallest coordinate arguments move towards one
	another. With three coordinates, their median gives the required
	control of both extremes and remains continuous when the
	ordering of the arguments changes.
	
	All coordinate indices are taken modulo \(3\).
	
	For \(\theta=(\theta_1,\theta_2,\theta_3)\in\mathbb R^3\),
	choose a permutation \(\sigma\) of \(\{1,2,3\}\) such that
	\(\theta_{\sigma(1)}\leq\theta_{\sigma(2)}\leq\theta_{\sigma(3)}\).
	Write
	\[
	\begin{aligned}
		a&=\min_i\theta_i=\theta_{\sigma(1)},\\
		m&=\med(\theta_1,\theta_2,\theta_3)=\theta_{\sigma(2)},\\
		b&=\max_i\theta_i=\theta_{\sigma(3)},
	\end{aligned}
	\qquad \Delta=b-a.
	\]
	The median \(m\) is the middle entry, with repeated values
	counted according to their multiplicity; equivalently,
	\(m=\theta_1+\theta_2+\theta_3-a-b\).
	The maximum and minimum are \(1\)-Lipschitz for the
	supremum norm. The identity for \(m\) therefore shows that
	\(a\), \(b\), \(m\) and \(\Delta\) are globally Lipschitz functions
	of \(\theta\). Finally, set
	\begin{equation}\label{eq:beta}
		\beta_i=d(\theta_{i+1}-m)-(\theta_i-m)
		=d\theta_{i+1}-\theta_i-(d-1)m
		\qquad(i=1,2,3).
	\end{equation}
	
	Let \(U_d\subset\PP^2\) be the set of
	points \([z]=[z_1:z_2:z_3]\) with \(z_i\ne0\) for all \(i\)
	for which there are real arguments \(\theta_i\) satisfying
	\(z_i=|z_i|e^{\mathrm i\theta_i}\) and
	\begin{equation}\label{eq:sector}
		\Delta=\max_i\theta_i-\min_i\theta_i<\frac{\pi}{2d}.
	\end{equation}
	Multiplication of \(z\) by a nonzero scalar adds the same real
	number to the three arguments, after a compatible choice of its
	argument. Neither their spread nor the quantities \(\beta_i\) in
	\eqref{eq:beta} change under such a translation. Thus the condition
	defining \(U_d\) depends only on the projective point. It is open:
	near a point with nonzero coordinates we may choose continuous local
	branches of all three arguments, and the strict inequality in
	\eqref{eq:sector} persists under a sufficiently small perturbation. Moreover,
	it contains the positive real locus
	\(\{[r_1:r_2:r_3]:r_i>0\}\).
	As in Corollary~\ref{cor:symmetry-count}, write
	\(A_{U_d}(q)=\{g(p_*):g\in G_d,\ L_q\cap g(U_d)\ne\varnothing\}\).
	The upper bound in \eqref{eq:sector} ensures positivity of
	the cosine factors in the radial equations; the choice of the
	median provides the corresponding control of the angular
	velocities. In what follows the arguments are real-valued
	variables, rather than angles taken modulo \(2\pi\). This ensures that the spread remains a continuous function
	along a trajectory. If two triples of arguments for the same
	representative satisfy \eqref{eq:sector}, their pairwise
	differences agree modulo \(2\pi\). Each such difference belongs to
	\(( -\delta,\delta )\), where
	\(\delta=\pi/(2d)\), and \(2\delta<2\pi\). Therefore, the
	two triples differ by a common integral multiple of \(2\pi\). Hence \(\Delta\)
	and each \(\beta_i\) have the same value for either choice; the
	phase \(e^{-\mathrm i(d-1)m}\) is unchanged as well. We fix an admissible triple at the initial point and continue
	it as a real-valued solution of the angular equations.
	No change of branch is required along the trajectory.
	
	\Needspace{5\baselineskip}
	\begin{lemma}\label{lem:median}
		With the notation above, \(\lvert\beta_i\rvert\leq d\Delta\)
		for every \(i\). Moreover, for every index \(i\) with
		\(\theta_i=b\) and every index \(j\) with \(\theta_j=a\),
		\begin{equation}\label{eq:extremes}
			\beta_i\leq -(b-m),\qquad
			\beta_j\geq m-a .
		\end{equation}
	\end{lemma}
	
	\begin{proof}
		Put \(u_k=\theta_k-m\). For distinct indices \(i,j\),
		the numbers \(u_i,u_j\) have opposite signs or one is zero:
		among three entries there is at most one strictly above the
		median and at most one strictly below it. Hence
		\[
		|\beta_i|=d|u_{i+1}|+|u_i|
		\leq d\bigl(|u_{i+1}|+|u_i|\bigr)
		=d|\theta_{i+1}-\theta_i|\leq d\Delta.
		\]
		This includes coincident arguments.
		
		The ordered list \((a,m,b)\) has at most one entry strictly
		above \(m\) and at most one entry strictly below \(m\), even
		when some entries coincide. Under our cyclic convention,
		\(k\) and \(k+1\) are distinct indices. If \(\theta_i=b\),
		then \(\theta_{i+1}>m\) would give
		\(\theta_i\geq\theta_{i+1}>m\), placing two entries strictly
		above \(m\). Thus \(\theta_{i+1}\leq m\). Likewise, if
		\(\theta_j=a\), then \(\theta_{j+1}<m\) would give
		\(\theta_j\leq\theta_{j+1}<m\), placing two entries strictly
		below \(m\). Thus \(\theta_{j+1}\geq m\). Since \(d\geq1\),
		substitution in \eqref{eq:beta} gives
		\begin{align*}
			\beta_i &= d(\theta_{i+1}-m)-(b-m)\leq -(b-m),\\
			\beta_j &= d(\theta_{j+1}-m)+(m-a)\geq m-a.\qedhere
		\end{align*}
	\end{proof}
	
	The restriction to three coordinates is used precisely in
	\eqref{eq:extremes}: the successor of an index attaining the
	maximum cannot lie strictly above the median, and the successor
	of an index attaining the minimum cannot lie strictly below it.
	
	\begin{remark}
		The only singular point in \(U_d\) is \(p_*\). Indeed, if
		\([z]\in U_d\) is singular, choose arguments as in
		\eqref{eq:sector}, write \(z_i=r_i e^{\mathrm i\theta_i}\)
		with \(r_i>0\), and take \(\lambda\) such that
		\(J_d(z)=\lambda z\); then \(\lambda\ne0\), since \(J_d(z)=0\)
		forces \(z=0\). Comparing the \(i\)-th components,
		\[
		e^{-\mathrm i(d-1)m}\lambda
		=\frac{r_{i+1}^d}{r_i}e^{\mathrm i\beta_i}
		\qquad(i=1,2,3).
		\]
		The three numbers \(\beta_i\) are arguments of the same nonzero
		complex number, hence congruent modulo \(2\pi\); since
		\(|\beta_i|\leq d\Delta<\pi/2\), they coincide.
		Their common value is at most \(-(b-m)\) and at least
		\(m-a\), by Lemma~\ref{lem:median}. Hence
		\[
		0\leq\Delta=(b-m)+(m-a)\leq0.
		\]
		Thus \(a=m=b\), and \(z_1/z_3\) and \(z_2/z_3\) are positive real
		numbers, and \eqref{eq:singular-set} forces \([z]=p_*\).
	\end{remark}
	
	We regard \((r,\theta)=(r_1,r_2,r_3,\theta_1,\theta_2,\theta_3)\)
	as real variables in the open set
	\(\Omega=(0,\infty)^3\times\mathbb R^3\).
	At each state put \(z_i=r_i e^{\mathrm i\theta_i}\), and use
	\(a,b,m,\Delta\) and \(\beta_i\) as above.
	We rotate \(J_d\) by \(e^{-\mathrm i(d-1)m}\). The choice of
	phase is dictated by \eqref{eq:beta}: for each \(i\),
	\[
	e^{-\mathrm i\theta_i}e^{-\mathrm i(d-1)m}z_{i+1}^{d}
	=r_{i+1}^{d}e^{\mathrm i\beta_i}.
	\]
	Multiplication by the nonzero phase preserves the projective
	tangent line of \(\J_d\); a radial term is annihilated
	by \(D\pi_z\). Writing
	\(\dot z=e^{-\mathrm i(d-1)m}J_d(z)-\Lambda z\), we have
	\[
	\dot r_i+\mathrm i r_i\dot\theta_i
	=r_{i+1}^d e^{\mathrm i\beta_i}-\Lambda r_i.
	\]
	We choose the real scalar \(\Lambda\) so that the positive
	normalisation \(r_1+r_2+r_3=1\) is preserved. Taking real parts
	and summing gives
	\[
	\sum_{i=1}^3\dot r_i
	=\sum_{i=1}^3r_{i+1}^d\cos\beta_i
	-\Lambda\sum_{i=1}^3r_i.
	\]
	Accordingly, set
	\[
	c_i=\cos\beta_i,\qquad
	\Lambda=\sum_{j=1}^3r_{j+1}^d c_j.
	\]
	The real and imaginary parts yield the system
	\begin{equation}\label{eq:polar}
		\dot r_i=r_{i+1}^d c_i-\Lambda r_i,\qquad
		\dot\theta_i=\frac{r_{i+1}^d}{r_i}\sin\beta_i
		\qquad (i=1,2,3).
	\end{equation}
	By construction, a solution \((r(t),\theta(t))\) of
	\eqref{eq:polar} yields a curve
	\(z(t)=(r_i(t)e^{\mathrm i\theta_i(t)})_{i=1}^3\) in
	\(\C^3\setminus\{0\}\) satisfying
	\begin{equation}\label{eq:lift}
		\dot z=e^{-\mathrm i(d-1)m}J_d(z)-\Lambda z .
	\end{equation}
	The system \eqref{eq:polar} is autonomous on \(\Omega\).
	The normalisation fixes the positive scale of \(z(t)\), but
	leaves its common argument free. This freedom is compatible with
	\eqref{eq:lift}. If \(z\) is replaced by \(e^{\mathrm i\varphi}z\),
	then \(m\) is replaced by \(m+\varphi\), whereas the moduli,
	the \(\beta_i\) and \(\Lambda\) remain unchanged. Homogeneity gives
	\[
	e^{-\mathrm i(d-1)(m+\varphi)}J_d(e^{\mathrm i\varphi}z)
	-\Lambda e^{\mathrm i\varphi}z
	=e^{\mathrm i\varphi}
	\bigl(e^{-\mathrm i(d-1)m}J_d(z)-\Lambda z\bigr).
	\]
	Thus the projected velocity is independent of the remaining
	choice of representative. Although the ordering of the arguments
	can change, the median remains Lipschitz where two arguments
	coincide. Consequently the right-hand side of \eqref{eq:polar}
	is locally Lipschitz on \(\Omega\), including these loci.
	Equation \eqref{eq:lift} also shows that the projected curve
	is tangent to \(\J_d\) wherever it is regular.
	
	Denote this projected real vector field by \(W_d\).
	At \(p_*\), the normalised moduli equal \(1/3\) and all
	\(\beta_i\) vanish. Equation~\eqref{eq:lift} then gives
	\(W_d(p_*)=0\).
	
	\Needspace{5\baselineskip}
	\begin{proposition}\label{prop:sector}
		The locally Lipschitz real field \(W_d\) on \(U_d\) has a
		positive flow defined for all time in \(U_d\). Every orbit
		starting at a regular point \(p\in U_d\) stays in \(L_p\)
		at every finite time and converges to \(p_*\). In particular,
		\(p_*\in\overline{L_p}\).
	\end{proposition}
	
	\Needspace{5\baselineskip}
	\begin{lemma}\label{lem:polar-existence}
		Suppose that \(r_i(0)>0\), \(\sum_i r_i(0)=1\), and
		\(\Delta(0)<\pi/(2d)\). The unique maximal solution of
		\eqref{eq:polar} with these initial data is defined for all
		\(t\geq0\), and \(\sum_i r_i(t)=1\). Moreover,
		\(a(t)\) is nondecreasing, \(b(t)\) is nonincreasing, and
		\[
		r_i(t)\geq r_i(0)e^{-t},\qquad
		c_i(t)\geq c_*:=\cos\bigl(d\Delta(0)\bigr)>0
		\qquad(t\geq0).
		\]
	\end{lemma}
	
	\begin{proof}
		As noted above, \(m\) and hence the \(\beta_i\) are globally
		Lipschitz functions of \(\theta\). The trigonometric functions
		are smooth, as are the powers of the \(r_i\) and the factors
		\(1/r_i\) on \(r_i>0\). Thus the right-hand side of
		\eqref{eq:polar} is locally Lipschitz on \(\Omega\), and the
		Picard--Lindel\"of theorem gives a unique maximal \(C^1\)
		solution with the specified initial data.
		
		Summing the equations for the \(r_i\) establishes the
		normalisation. If \(S(t)=\sum_i r_i(t)\), then
		\(\dot S=\Lambda(1-S)\) and \(S(0)=1\). Since \(\Lambda(t)\)
		is continuous along the solution, integration gives
		\[
		S(t)-1=(S(0)-1)\exp\!\left(-\int_0^t\Lambda(s)\,ds\right)=0 .
		\]
		Thus \(S(t)=1\) throughout the interval of existence; the radial
		normalisation does not rely on an a priori bound on the lifetime
		of the solution.
		
		To prove positive invariance of the sector, write
		\(a(t),b(t),\Delta(t)\) for the extremal values and the spread
		of \(\theta(t)\). The functions \(a,b\) are locally Lipschitz
		because the \(\theta_i\) are \(C^1\), hence differentiable
		almost everywhere. If \(b\) is differentiable at \(t\) and
		\(\theta_i(t)=b(t)\), then \(b-\theta_i\) has a local
		minimum equal to zero, whence
		\(\dot b(t)=\dot\theta_i(t)\). Similarly,
		\(\dot a(t)=\dot\theta_j(t)\) for every index satisfying
		\(\theta_j(t)=a(t)\) at a time when \(a\) is differentiable.
		These identities also cover ties between the arguments.
		
		As long as \(d\Delta<\pi/2\), the bound on \(\beta_i\)
		in Lemma~\ref{lem:median} gives \(|\beta_i|<\pi/2\), so
		\(\sin\beta_i\) has the sign of \(\beta_i\) and
		\(c_i=\cos\beta_i>0\). By
		\eqref{eq:extremes} and \eqref{eq:polar}, the derivative
		\(\dot\theta_i\) is nonpositive at any index attaining the
		maximum, and \(\dot\theta_j\) is nonnegative at any index
		attaining the minimum. Thus, as long as the sector condition holds,
		\[
		\dot b(t)\leq0,\qquad \dot a(t)\geq0
		\qquad\text{for almost every }t.
		\]
		Local absolute continuity now shows that \(b\) is nonincreasing
		and \(a\) is nondecreasing on every such interval. To verify positive invariance of the sector, suppose that
		\(t_1\) is the first forward time at which
		\(d\Delta(t_1)=\pi/2\). Then
		\(\Delta\) would be nonincreasing on \([0,t_1)\); continuity
		would imply
		\[
		\frac{\pi}{2d}=\Delta(t_1)
		\leq\Delta(0)<\frac{\pi}{2d},
		\]
		a contradiction. Therefore, throughout the forward interval of
		existence,
		\[
		a(0)\leq\theta_i(t)\leq b(0),\qquad
		c_i(t)\geq c_*:=\cos\bigl(d\Delta(0)\bigr)>0.
		\]
		Here the cosine bound uses \(|\beta_i(t)|\leq d\Delta(t)\leq d\Delta(0)\)
		and the decrease of cosine on \([0,\pi/2]\).
		
		We establish global existence for \eqref{eq:polar} by excluding
		a finite right endpoint of the maximal interval. The normalisation
		gives \(0<r_i\leq1\), hence
		\(\sum_i r_i^d\leq\sum_i r_i=1\), and \(0<c_j\leq1\) by the
		previous paragraph. Thus
		\[
		0<\Lambda=\sum_{j=1}^3r_{j+1}^dc_j
		\leq\sum_{j=1}^3r_{j+1}^d
		\leq\sum_{j=1}^3r_{j+1}=1,
		\]
		and \eqref{eq:polar} together with \(c_i\geq c_*\) yields
		\begin{equation}\label{eq:modulus-comparison}
			\dot r_i\geq c_*r_{i+1}^d-r_i\geq-r_i.
		\end{equation}
		In particular,
		\[
		\frac{d}{dt}(e^t r_i(t))
		=e^t(\dot r_i(t)+r_i(t))
		\geq e^t c_*r_{i+1}(t)^d\geq0,
		\]
		so that
		\[
		r_i(t)\geq r_i(0)e^{-t}>0.
		\]
		This estimate is not uniform as \(t\to\infty\), but it is
		sufficient to exclude loss of positivity in finite time.
		If the maximal forward time
		were \(T<\infty\), each \(r_i(t)\) would belong to
		\([r_i(0)e^{-T},1]\) for \(0\leq t<T\), and each
		\(\theta_i(t)\) to \([a(0),b(0)]\). These bounds keep the
		state in a compact subset of \(\Omega\), on which the
		velocity is bounded. The state is therefore Lipschitz on
		\([0,T)\) and has a limit as \(t\uparrow T\). Local existence
		at that limit extends the solution past \(T\). We therefore
		have a solution for every \(t\geq0\).
	\end{proof}
	
	\Needspace{5\baselineskip}
	\begin{lemma}\label{lem:angular-convergence}
		For the solution in Lemma~\ref{lem:polar-existence}, there
		are constants \(\eta,\kappa>0\), depending only on \(d\)
		and \(\Delta(0)\), such that
		\[
		r_i(t)\geq\eta\quad(i=1,2,3),\qquad
		\Delta(t)\leq\Delta(2)e^{-\kappa(t-2)}
		\qquad(t\geq2).
		\]
		One may take \(\kappa=2\eta^d/\pi\). Furthermore,
		\(\lvert c_i(t)-1\rvert\leq|\beta_i(t)|\leq d\Delta(t)\) for all \(t\geq0\).
	\end{lemma}
	
	\begin{proof}
		A positive lower bound for the moduli, uniform in time, makes
		the angular inequalities quantitative. Such a bound follows
		from the cyclic form of the radial equations. The inequality from
		\eqref{eq:modulus-comparison},
		\[
		\dot r_i+r_i\geq c_*r_{i+1}^d,
		\]
		transfers a positive lower bound to the preceding coordinate.
		At time \(t-2\), some modulus is at least \(1/3\). We propagate
		this bound through the other two coordinates over successive
		unit intervals. The initial index may depend on \(t\); the constants below
		are independent of that choice.
		
		Fix \(t\geq2\) and put \(s=t-2\). Since the moduli sum to
		one, some \(r_j(s)\) is at least \(1/3\). The inequality
		\(\dot r_j\geq-r_j\) gives
		\[
		r_j(\tau)\geq \frac{e^{-(\tau-s)}}{3}
		\geq A_0:=\frac{e^{-2}}3,
		\qquad s\leq\tau\leq t.
		\]
		The cyclic equation for \(r_{j-1}\) contains \(r_j^d\).
		The lower bound just obtained for \(r_j\) holds throughout
		\([s,t]\), and therefore in particular on its first unit
		subinterval. Consequently,
		\[
		\dot r_{j-1}+r_{j-1}\geq c_*A_0^d
		\qquad\text{on }[s,t].
		\]
		Multiplying by \(e^\tau\) and integrating from \(s\) to
		\(\tau\in[s,t]\), we obtain
		\[
		\begin{aligned}
			r_{j-1}(\tau)
			&\geq e^{-(\tau-s)}r_{j-1}(s)
			+c_*A_0^d\int_s^\tau e^{-(\tau-v)}\,dv\\
			&\geq c_*A_0^d(1-e^{-(\tau-s)}).
		\end{aligned}
		\]
		In particular, for \(s+1\leq\tau\leq t\),
		\[
		r_{j-1}(\tau)\geq
		A_1:=c_*A_0^d(1-e^{-1}).
		\]
		Since the bound for \(r_{j-1}\) holds throughout
		\([s+1,t]\), it can be inserted into the equation for
		\(r_{j-2}\) on the entire second unit interval. Thus
		\[
		\dot r_{j-2}+r_{j-2}\geq c_*A_1^d.
		\]
		Integrating from \(s+1\) to \(t=s+2\) gives
		\[
		\begin{aligned}
			r_{j-2}(t)
			&\geq e^{-1}r_{j-2}(s+1)
			+c_*A_1^d\int_{s+1}^t e^{-(t-v)}\,dv\\
			&\geq A_2:=c_*A_1^d(1-e^{-1}).
		\end{aligned}
		\]
		At time \(t\) these estimates cover all three coordinates. Thus
		\begin{equation}\label{eq:persistence}
			r_i(t)\geq\eta:=\min\{A_0,A_1,A_2\}>0,
			\qquad t\geq2,\quad i=1,2,3.
		\end{equation}
		The constants can also be written explicitly. Set
		\(B=c_*(1-e^{-1})\). Since \(0<B<1\) and \(0<A_0<1\),
		\[
		A_1=BA_0^d<A_0,\qquad A_2=BA_1^d<A_1.
		\]
		Consequently,
		\[
		0<A_2<A_1<A_0,\qquad
		\eta=A_2=B^{d+1}A_0^{d^2},\qquad
		\kappa=\frac{2\eta^d}{\pi}
		=\frac2\pi B^{d(d+1)}A_0^{d^3}.
		\]
		In particular, \(\kappa>0\). Neither \(\eta\) nor
		\(\kappa\) depends on the initial moduli: the initial
		normalisation alone supplied the coordinate of size at least
		\(1/3\). They do depend on the angular spread through \(c_*\),
		so these constants are not uniform on all of \(U_d\).
		
		The estimate \eqref{eq:persistence} is what turns the sign
		control of Lemma~\ref{lem:median} into a quantitative contraction:
		for \(t\geq2\), the bounds \(\eta\leq r_i\leq1\) give
		\[
		\frac{r_{i+1}(t)^d}{r_i(t)}\geq \eta^d.
		\]
		Thus the coefficients of the sine terms in the angular equations
		are uniformly bounded below after time two. We also use the elementary
		inequality
		\[
		\sin u\geq\frac{2u}{\pi},\qquad 0\leq u\leq\frac\pi2,
		\]
		which follows from concavity. At an index \(i\) where \(\theta_i=b\), we have
		\(0\leq b-m\leq-\beta_i<\pi/2\); at an index \(j\) where
		\(\theta_j=a\), we have
		\(0\leq m-a\leq\beta_j<\pi/2\).
		Since sine is increasing on \([0,\pi/2]\),
		\eqref{eq:polar}, \eqref{eq:persistence} and the
		sine inequality give
		\[
		\begin{aligned}
			\dot\theta_i
			&=\frac{r_{i+1}^d}{r_i}\sin\beta_i
			\leq-\frac{r_{i+1}^d}{r_i}\sin(b-m)
			\leq-\kappa(b-m),\\
			\dot\theta_j
			&=\frac{r_{j+1}^d}{r_j}\sin\beta_j
			\geq\frac{r_{j+1}^d}{r_j}\sin(m-a)
			\geq\kappa(m-a).
		\end{aligned}
		\]
		The maximum and minimum of finitely many \(C^1\) functions
		are locally absolutely continuous. At almost every time when
		they are differentiable, their derivatives agree with the
		derivatives of the active coordinates, even when several
		coordinates share the extremum, as shown in
		Lemma~\ref{lem:polar-existence}.
		Thus
		\[
		\dot b(t)\leq-\kappa(b(t)-m(t)),\qquad
		\dot a(t)\geq\kappa(m(t)-a(t))
		\]
		for almost every \(t\geq2\). Subtracting gives
		\[
		\dot\Delta(t)\leq-\kappa\bigl((b-m)+(m-a)\bigr)
		=-\kappa\Delta(t)
		\qquad\text{a.e. on }[2,\infty).
		\]
		Equivalently,
		\[
		\frac{d}{dt}\bigl(e^{\kappa(t-2)}\Delta(t)\bigr)\leq0
		\qquad\text{for almost every }t\geq2.
		\]
		Absolute continuity and integration yield
		\begin{equation}\label{eq:phase-decay}
			\Delta(t)\leq\Delta(2)e^{-\kappa(t-2)},\qquad t\geq2.
		\end{equation}
		Finally, the mean value theorem and Lemma~\ref{lem:median}
		show that, for every \(t\geq0\),
		\[
		|c_i(t)-1|\leq|\beta_i(t)|\leq d\Delta(t).
		\]
		Moreover, since \(c_i=\cos\beta_i\), the inequality
		\(1-\cos u\leq u^2/2\) gives
		\begin{equation}\label{eq:quadratic-cosine}
			|c_i(t)-1|=1-\cos\beta_i(t)
			\leq\frac{\beta_i(t)^2}{2}
			\leq\frac{d^2}{2}\Delta(t)^2.
		\end{equation}
	\end{proof}
	
	The restriction of the affine field \(X_d\) in
	\eqref{eq:affine-generator} to \((0,\infty)^2\) is a real
	vector field with an equilibrium at \((1,1)\). For a real-valued
	\(C^1\) function \(V\) on this quadrant, we use the standard
	notation \(\mathcal L_{X_d}V=dV(X_d)\) for its derivative along
	the restricted field. We choose \(V(x,y)=f(x)+g(y)\), where
	\[
	f'(x)=x^d-1,\qquad yg'(y)=y^d-1.
	\]
	With this choice the mixed terms in \(\mathcal L_{X_d}V\)
	cancel. The normalisations \(f(1)=g(1)=0\) determine the
	functions in the following lemma.
	
	\Needspace{5\baselineskip}
	\begin{lemma}\label{lem:lyapunov}
		For \(x,y>0\), define \(V(x,y)=f(x)+g(y)\), where
		\[
		f(x)=\frac{x^{d+1}-1}{d+1}-(x-1),\qquad
		g(y)=\frac{y^d-1}{d}-\log y .
		\]
		Then \(V\geq0\), with equality only at \((1,1)\), and
		\[
		(\mathcal L_{X_d}V)(x,y)\leq-V(x,y)
		\qquad(x,y>0).
		\]
	\end{lemma}
	
	\begin{proof}
		Both summands and their first derivatives vanish at \(1\):
		\[
		\begin{gathered}
			f(1)=g(1)=f'(1)=g'(1)=0,\\
			f'(x)=x^d-1,\qquad g'(y)=\frac{y^d-1}{y},\\
			f''(x)=dx^{d-1}>0,\qquad
			g''(y)=(d-1)y^{d-2}+y^{-2}>0.
		\end{gathered}
		\]
		Both \(f\) and \(g\) are therefore strictly convex on
		\((0,\infty)\), and \(V\geq0\), with equality only at
		\((1,1)\).
		
		To calculate the derivative along \(X_d\), we use
		\(yg'(y)=y^d-1\) and the identities
		\[
		\begin{aligned}
			y^d-x^{d+1}&=(y^d-1)-(x^{d+1}-1),\\
			1-x^dy&=-(y-1)-y(x^d-1).
		\end{aligned}
		\]
		It follows that
		\begin{align}
			(\mathcal L_{X_d}V)(x,y)
			&= (x^d-1)(y^d-x^{d+1})+g'(y)(1-x^dy) \notag\\
			&= -(x^d-1)(x^{d+1}-1)-(y-1)g'(y).
			\label{eq:lyapunov-identity}
		\end{align}
		
		We next compare the two terms on the right of
		\eqref{eq:lyapunov-identity} with \(f\) and \(g\).
		For a differentiable convex function \(h\), the supporting-line
		inequality at \(u>0\), evaluated at the point \(1\), gives
		\[
		h(1)\geq h(u)+h'(u)(1-u).
		\]
		Applying this inequality to \(f\) and \(g\), both of which
		vanish at \(1\), gives
		\[
		f(x)\leq(x-1)f'(x),\qquad
		g(y)\leq(y-1)g'(y).
		\]
		Furthermore,
		\(x^{d+1}-1=(x-1)(1+x+\cdots+x^d)\). The factor
		\(1+x+\cdots+x^d\) is at least \(1\), while
		\((x-1)(x^d-1)\geq0\). Hence
		\[
		(x^d-1)(x^{d+1}-1)
		\geq(x-1)f'(x)\geq f(x).
		\]
		The convexity inequality for \(g\) also gives
		\[
		(y-1)g'(y)\geq g(y).
		\]
		Together with the preceding estimate and
		\eqref{eq:lyapunov-identity}, we obtain
		\begin{equation}\label{eq:lyapunov}
			(\mathcal L_{X_d}V)(x,y)\leq -V(x,y)
			\qquad(x,y>0),
		\end{equation}
		as required.
	\end{proof}
	
	\begin{proof}[Proof of Proposition~\ref{prop:sector}]
		Choose a representative \(z^0=(z_1^0,z_2^0,z_3^0)\) of \(p\)
		and real arguments \(\theta_i(0)\) of \(z_i^0\) satisfying
		\eqref{eq:sector}. Normalise the representative by setting
		\[
		z(0)=\frac{z^0}{\sum_{k=1}^3 |z_k^0|},\qquad
		r_i(0)=|z_i(0)|.
		\]
		This positive rescaling preserves both the projective point
		and the chosen real arguments. Hence \(\pi(z(0))=p\), the initial
		moduli are positive, and \(\sum_i r_i(0)=1\).
		
		Let \((r(t),\theta(t))\) be the global solution of
		\eqref{eq:polar} supplied by Lemma~\ref{lem:polar-existence}.
		Set
		\[
		z(t)=(r_i(t)e^{\mathrm i\theta_i(t)})_{i=1}^3,
		\qquad \gamma(t)=[z(t)].
		\]
		Then \(z(t)\) satisfies \eqref{eq:lift}, and its projectivisation
		is the positive orbit of \(W_d\) through \(p\).
		
		The curve can be studied throughout \([0,\infty)\) in the
		affine coordinates \(w=(w_1,w_2)=(z_1/z_3,z_2/z_3)\), since
		\(r_3>0\). These coordinates also identify the leaf containing
		the curve.
		Differentiating \(w_k=z_k/z_3\) gives
		\(\dot w_k=(\dot z_k-w_k\dot z_3)/z_3\), so the radial
		term in \eqref{eq:lift} cancels. Using the components
		\(J_d(z)=(z_2^d,z_3^d,z_1^d)\), we obtain
		\[
		\dot w(t)=\alpha(t)X_d(w(t)),\qquad
		\alpha(t)=e^{-\mathrm i(d-1)m(t)}z_3(t)^{d-1},
		\]
		with \(\alpha\) continuous and nowhere zero, because
		\(r_3(t)>0\). Although the phase defining \(\alpha\) depends
		only Lipschitz-continuously on the arguments, the resulting
		equation has the uniqueness property established in
		Lemma~\ref{lem:tangent-curve}(2). Since \(w(0)\) is regular,
		that lemma shows that \(\gamma(t)\) is regular for every finite
		\(t\geq0\). Moreover, \(\dot\gamma(t)\) lies in the complex
		tangent line of \(\J_d\), and Lemma~\ref{lem:tangent-curve}(1)
		then yields
		\[
		\gamma([0,\infty))\subset L_p.
		\]
		Thus the construction gives a curve in the original regular
		leaf for every finite time. Any singular limit is approached
		as \(t\to\infty\); to identify it, we estimate the moduli
		and the differences of arguments.
		
		By Lemma~\ref{lem:angular-convergence}, and in particular
		\eqref{eq:persistence}, there are constants \(\eta,\kappa>0\)
		such that all three moduli stay above \(\eta\) for
		\(t\geq2\) and
		\[
		\Delta(t)\leq\Delta(2)e^{-\kappa(t-2)}.
		\]
		Thus the arguments become asymptotically equal. It remains to
		show that the moduli do likewise.
		
		Put \(x=r_1/r_3=|w_1|\) and \(y=r_2/r_3=|w_2|\).
		Differentiating and using \eqref{eq:polar}, the radial terms
		involving \(\Lambda\) cancel:
		\[
		\begin{aligned}
			\dot x
			&=\frac{\dot r_1-x\dot r_3}{r_3}
			=\frac{c_1r_2^d-xc_3r_1^d}{r_3},\\
			\dot y
			&=\frac{\dot r_2-y\dot r_3}{r_3}
			=\frac{c_2r_3^d-yc_3r_1^d}{r_3}.
		\end{aligned}
		\]
		Using \(r_1=xr_3\) and \(r_2=yr_3\), we obtain
		\begin{equation}\label{eq:ratios}
			\dot x=r_3^{d-1}(c_1y^d-c_3x^{d+1}),\qquad
			\dot y=r_3^{d-1}(c_2-c_3x^dy).
		\end{equation}
		The cancellation of the normalising term \(\Lambda\) is the
		reason for using the ratios \(x=r_1/r_3\), \(y=r_2/r_3\).
		When \(c_1=c_2=c_3=1\), equation \eqref{eq:ratios} is
		\(\dot q=r_3^{d-1}X_d(q)\) on the positive real quadrant.
		In particular, if
		\(\Delta(0)=0\), then the angular equations vanish, the
		coordinates retain a common argument, and the ratios follow
		this positive real flow after the change of time
		\(ds/dt=r_3^{d-1}\).
		
		For general initial arguments the quantities \(c_i-1\) do
		not vanish, but Lemma~\ref{lem:angular-convergence} shows
		that they tend to zero exponentially. We use
		Lemma~\ref{lem:lyapunov} to control the resulting
		perturbation of the positive real system.
		
		Put \(q(t)=(x(t),y(t))\). For \(t\geq2\), the bounds
		\(\eta\leq r_i(t)\leq1\) imply that
		\[
		q(t)\in\mathcal R=[\eta,\eta^{-1}]^2.
		\]
		Equations \eqref{eq:ratios} can be written as
		\[
		\begin{aligned}
			\dot q(t)&=r_3(t)^{d-1}X_d(q(t))+E(t),\\
			E(t)&=r_3(t)^{d-1}\bigl(
			(c_1-1)y^d-(c_3-1)x^{d+1},\,
			(c_2-1)-(c_3-1)x^dy\bigr),
		\end{aligned}
		\]
		where \(x=x(t)\), \(y=y(t)\) and \(c_i=c_i(t)\) on the
		right-hand side. Thus \(E(t)\) is a real vector along
		\(q(t)\), and its contribution to the derivative of \(V\)
		is given by the differential
		\[
		dV_{(x,y)}(\xi)=f'(x)\xi_1+g'(y)\xi_2
		\qquad(\xi=(\xi_1,\xi_2)\in\mathbb R^2).
		\]
		Consequently,
		\[
		\begin{aligned}
			|dV_{q(t)}(E(t))|
			&\leq r_3(t)^{d-1}\max_i|c_i(t)-1|\,
			\bigl[|f'(x)|(y^d+x^{d+1})\\
			&\hspace{91pt}{}+|g'(y)|(1+x^dy)\bigr].
		\end{aligned}
		\]
		The bracket is bounded on the compact rectangle
		\(\mathcal R\), say by \(M_{d,\eta}\), uniformly in time.
		Since
		\(r_3(t)^{d-1}\leq1\),
		\eqref{eq:quadratic-cosine} gives
		\[
		|dV_{q(t)}(E(t))|
		\leq \frac{d^2}{2}M_{d,\eta}\Delta(t)^2
		=C\Delta(t)^2,\qquad t\geq2,
		\]
		where \(C\) depends only on \(d\) and \(\eta\).
		Set \(\mu=\eta^{d-1}>0\). The same persistence estimate
		that bounds the error term also gives
		\[
		\mu\leq r_3(t)^{d-1}\leq1 \qquad(t\geq2).
		\]
		The lower bound and the decaying error have been obtained
		without assuming convergence of \(q(t)\). Applying the
		chain rule to the Lyapunov function gives
		\[
		\begin{aligned}
			\frac{d}{dt}V(q(t))
			&=dV_{q(t)}(\dot q(t))\\
			&=r_3(t)^{d-1}(\mathcal L_{X_d}V)(q(t))
			+dV_{q(t)}(E(t)).
		\end{aligned}
		\]
		Since \(r_3(t)^{d-1}\geq\mu\), the inequality
		\eqref{eq:lyapunov} and \(V\geq0\) yield
		\begin{equation}\label{eq:perturbed-lyapunov}
			\frac{d}{dt}V(q(t))
			\leq-\mu V(q(t))+C\Delta(t)^2,\qquad t\geq2.
		\end{equation}
		Write \(W(t)=V(q(t))\). Integrating
		\eqref{eq:perturbed-lyapunov} and using
		\eqref{eq:phase-decay}, we obtain
		\[
		\begin{aligned}
			0\leq W(t)
			&\leq e^{-\mu(t-2)}W(2)\\
			&\quad+C\Delta(2)^2\int_2^t
			e^{-\mu(t-\tau)-2\kappa(\tau-2)}\,d\tau.
		\end{aligned}
		\]
		For \(2\leq\tau\leq t\), the nonnegative numbers
		\(t-\tau\) and \(\tau-2\) sum to \(t-2\). Hence
		\[
		\mu(t-\tau)+2\kappa(\tau-2)
		\geq\min\{\mu,2\kappa\}(t-2).
		\]
		This estimate also covers \(\mu=2\kappa\), without requiring
		division by the difference of the two exponents. Hence
		\[
		\begin{aligned}
			0\leq W(t)
			&\leq e^{-\mu(t-2)}W(2)\\
			&\quad+C\Delta(2)^2(t-2)
			e^{-\min\{\mu,2\kappa\}(t-2)}\longrightarrow0.
		\end{aligned}
		\]
		Thus \(V(x(t),y(t))\to0\). To deduce convergence of the
		ratios, note that for any \(\varepsilon>0\) for which
		\(\mathcal R\setminus B_\varepsilon(1,1)\) is nonempty,
		compactness and strict positivity give
		\[
		\min_{(x,y)\in\mathcal R\setminus B_\varepsilon(1,1)}
		V(x,y)>0.
		\]
		Since \(V(x(t),y(t))\to0\), the ratios eventually belong to
		\(B_\varepsilon(1,1)\). We conclude that
		\[
		x(t)\longrightarrow1,\qquad y(t)\longrightarrow1.
		\]
		
		It remains to pass from the moduli and arguments to projective
		coordinates. The two differences \(\theta_i-\theta_3\) for \(i=1,2\)
		are bounded in absolute value by \(\Delta(t)\), hence
		tend to zero. In the affine chart,
		\[
		\gamma(t)=
		[x(t)e^{\mathrm i(\theta_1(t)-\theta_3(t))}:
		y(t)e^{\mathrm i(\theta_2(t)-\theta_3(t))}:1]
		\longrightarrow [1:1:1]=p_* .
		\]
		The curve stays in \(L_p\) for every finite time, while its
		projective limit is the singular point \(p_*\). Consequently,
		\(p_*\in\overline{L_p}\), as asserted.
	\end{proof}
	
	\Needspace{5\baselineskip}
	\begin{corollary}\label{cor:sector-rate}
		Let \(p\in U_d\setminus\Sing(\J_d)\), and let \(\gamma_p(t)\),
		\(t\geq0\), be the projective curve of \(W_d\) through \(p\). There are constants \(C_p,\sigma_p>0\)
		such that
		\[
		\operatorname{dist}_{FS}(\gamma_p(t),p_*)
		\leq C_p e^{-\sigma_p t}\qquad(t\geq0),
		\]
		where \(\operatorname{dist}_{FS}\) is the Fubini--Study distance
		and \(t\) is the positive-flow time of the normalised field
		\(W_d\).
	\end{corollary}
	
	\begin{proof}
		Retain the notation of the proof of Proposition~\ref{prop:sector}.
		Put \(s=t-2\) for \(t\geq2\), and set
		\(a_0=\min\{\mu,2\kappa\}>0\). The integral estimate obtained
		there gives
		\[
		0\leq W(t)\leq e^{-\mu s}W(2)
		+C\Delta(2)^2 s e^{-a_0s}
		\leq C_0 e^{-a_0s/2}
		\]
		for a suitable constant \(C_0>0\). Both \(x(t)\) and \(y(t)\)
		remain in the compact interval \([\eta,\eta^{-1}]\).
		Since \(f''\) and \(g''\) are strictly positive there,
		there is a number \(m_0>0\) such that
		\[
		f''(x)\geq m_0,\qquad g''(y)\geq m_0
		\qquad (x,y\in[\eta,\eta^{-1}]).
		\]
		Taylor's theorem at \(1\), using
		\(f(1)=g(1)=f'(1)=g'(1)=0\), therefore gives
		\[
		V(x,y)\geq\frac{m_0}{2}
		\bigl((x-1)^2+(y-1)^2\bigr)
		\qquad((x,y)\in\mathcal R).
		\]
		It follows that
		\[
		|x(t)-1|+|y(t)-1|\leq C_1e^{-a_0s/4}.
		\]
		The differences of arguments are bounded by
		\(\Delta(t)\leq\Delta(2)e^{-\kappa s}\). In the affine
		chart \(z_3\ne0\),
		\[
		w(t)=\bigl(x(t)e^{\mathrm i(\theta_1-\theta_3)},
		y(t)e^{\mathrm i(\theta_2-\theta_3)}\bigr).
		\]
		Using \(\lvert e^{\mathrm iu}-1\rvert\leq\lvert u\rvert\)
		and the boundedness of \(x\) and \(y\), we obtain
		\[
		|w_1(t)-1|+|w_2(t)-1|
		\leq C_2e^{-a_0s/4}.
		\]
		On the compact set containing \(w(t)\) for \(t\geq2\),
		the Fubini--Study distance to \(p_*\) is bounded by a
		constant multiple of the Euclidean distance to \((1,1)\).
		We obtain the asserted estimate for \(t\geq2\), with
		\(\sigma_p=a_0/4\). Continuity of \(\gamma_p\) on
		\([0,2]\) allows us to enlarge \(C_p\) and include the
		initial interval.
	\end{proof}
	
	The estimates also give a Lyapunov function in the form of
	Proposition~\ref{prop:lyapunov-criterion}. With the notation
	of the proof, choose \(A\geq C/\kappa\) and put
	\[
	H=V(r_1/r_3,r_2/r_3)+A\Delta^2.
	\]
	On the compact part of the sector reached after time two,
	\eqref{eq:perturbed-lyapunov} and the angular inequality give
	\[
	\frac{dH}{dt}
	\leq-\mu V+(C-2A\kappa)\Delta^2
	\leq-\min\{\mu,\kappa\}H
	\]
	almost everywhere. This compact set is contained in \(U_d\),
	since all moduli are at least \(\eta\) and the angular spread
	is at most its initial value, strictly below the sector width.
	The function \(H\) is locally Lipschitz and vanishes there
	only at \(p_*\). Thus the radial and angular estimates verify
	the convergence condition in the general surface criterion.
	
	\section{The coordinate lines and the singularity count}
	
	We now verify the finite jet alternative of
	Theorem~\ref{thm:propagation-count}. Let \(H=\{z_2=0\}\).
	The local calculation below gives the leading coefficients in
	Lemma~\ref{lem:jet-entry}; the diagonal symmetries select the
	admissible angular regions. Counting their singular points
	then gives the required finite sets.
	
	\Needspace{5\baselineskip}
	\begin{lemma}\label{lem:coordinate-line}
		Every point of the projective line
		\(H=\{[z]\in\PP^2 : z_2=0\}\) is regular for \(\J_d\).
		For each \(p\in H\), the closure of the leaf \(L_p\) through
		\(p\) meets \(\Sing(\J_d)\).
	\end{lemma}
	
	\begin{proof}
		If \([z]\in H\) were singular, then \(J_d(z)=\lambda z\)
		for some \(\lambda\in\C\). Since \(z_2=0\), the second
		component gives \(z_3^d=\lambda z_2=0\). The third then
		gives \(z_1^d=\lambda z_3=0\). This would make \(z=0\),
		which is impossible for a projective representative.
		
		The point \([1:0:0]\) is the only point of \(H\) outside the
		chart \(z_3\ne0\). The cyclic symmetry sends it to
		\([0:0:1]\in H\) and carries its leaf to the leaf through
		\([0:0:1]\). We may therefore work in the affine chart and
		consider points \([u:0:1]\), \(u\in\C\).
		The map induced by \(D_\zeta\) sends \([u:0:1]\) to
		\([\zeta u:0:1]\). If \(u\ne0\), take a
		real argument \(\psi\) of \(u\) and choose an integer \(k\)
		nearest to \(-N\psi/(2\pi)\). With \(\zeta=e^{2\pi\mathrm i k/N}\),
		the number \(\zeta u\) has argument \(\psi+2\pi k/N\),
		whose absolute value is at most \(\pi/N\). After applying
		this symmetry, either \(u=0\) or we may choose an argument
		\(\vartheta\) of \(u\) satisfying
		\begin{equation}\label{eq:angular-margin}
			|\vartheta|\leq\frac{\pi}{N}<\frac{\pi}{2d}.
		\end{equation}
		Indeed, \(N-2d=d^2-d+1>0\) for \(d\geq2\).
		(The same reduction, with the additional use of complex
		conjugation, is the first step of \cite{CF}.)
		
		From such a point, we reach \(U_d\) by following the projection
		of an integral curve of \(J_d\). Let \(z(t)\) be the local
		solution of \(\dot z=J_d(z)\)
		with \(z(0)=(u,0,1)\) and real time \(t\).
		The initial equations give
		\(\dot z_1(0)=0\), \(\dot z_2(0)=1\) and
		\(\dot z_3(0)=u^d\). Moreover,
		\(\dot z_1=z_2^d\) and \(z_2(t)=t+O(t^2)\), so
		\[
		\dot z_1(t)=t^d+O(t^{d+1}),\qquad
		z_1(t)-u=\int_0^t\bigl(s^d+O(s^{d+1})\bigr)\,ds.
		\]
		Together with the first derivatives of \(z_2\) and
		\(z_3\), this gives the expansions
		\[
		\begin{aligned}
			z_1(t)&=u+\frac{t^{d+1}}{d+1}+O(t^{d+2}),\\
			z_2(t)&=t+O(t^2),\\
			z_3(t)&=1+u^dt+O(t^2)
		\end{aligned}
		\qquad(t\downarrow0).
		\]
		In particular,
		\[
		z_1(t)\longrightarrow u,\qquad
		\frac{z_2(t)}t\longrightarrow1,\qquad
		z_3(t)\longrightarrow1.
		\]
		
		If \(u\ne0\), the limits of \(z_1(t)\), \(z_2(t)/t\), and
		\(z_3(t)\) are all nonzero. For sufficiently small \(t>0\)
		these three quantities admit continuous choices of arguments
		tending to \(\vartheta,0,0\), respectively, by taking a
		local branch of the argument near each limit. The argument
		chosen for \(z_2(t)/t\) is also an argument of \(z_2(t)\),
		since \(t>0\). The spread of these arguments tends to
		\(\max\{\vartheta,0\}-\min\{\vartheta,0\}=|\vartheta|\).
		By \eqref{eq:angular-margin}, the spread is smaller than
		\(\pi/(2d)\) for sufficiently small \(t>0\).
		Thus \([z(t)]\in U_d\).
		
		If \(u=0\), uniqueness and the reality of \(J_d\)
		show that \(z(t)\) is real. The preceding expansions become
		\[
		z_1(t)=\frac{t^{d+1}}{d+1}+O(t^{d+2}),\qquad
		z_2(t)=t+O(t^2),\qquad z_3(t)=1+O(t^2).
		\]
		All three coordinates are therefore positive for sufficiently
		small \(t>0\); their arguments can be chosen equal to zero,
		and \([z(t)]\in U_d\).
		
		In both cases, the projective trajectory \([z(t)]\) starts at
		the regular point \([u:0:1]\) and stays regular for small
		time, by openness of the regular locus. Its velocity is
		\(D\pi_{z(t)}J_d(z(t))\), tangent to \(\J_d\). It follows
		from Lemma~\ref{lem:tangent-curve}(1) that the points reached
		in \(U_d\) belong to the same leaf as \([u:0:1]\).
		Proposition~\ref{prop:sector} then places \(p_*\) in the
		closure of that leaf. Finally, applying the inverse of the
		projective symmetry used at the start gives a singular point
		in the closure of the leaf through the original point of \(H\).
	\end{proof}
	
	In the affine chart, the functions in
	Lemma~\ref{lem:jet-entry} are \(f_1=w_1\), \(f_2=w_2\),
	and \(f_3=1\). The set \(\Theta\) consists of triples of
	arguments with spread less than \(\pi/(2d)\). At
	\([u:0:1]\), \(u\ne0\), their orders and leading
	coefficients along the local flow of \(X_d\) are
	\((0,1,0)\) and \((u,1,1)\). At \(u=0\) they are
	\((d+1,1,0)\) and \((1/(d+1),1,1)\), since
	\(\dot w_2=1-w_1^dw_2\) and
	\(\dot w_1=w_2^d-w_1^{d+1}\).
	Thus the preceding entry calculation is precisely a finite jet
	condition, including the higher contact at the coordinate point.
	The median is used to establish attraction inside \(U_d\),
	whereas these leading coefficients establish entry into it.
	
	The proof of Lemma~\ref{lem:coordinate-line} applies to
	every diagonal symmetry satisfying the strict angular
	condition. Counting these symmetries gives a lower bound
	independent of the chosen point of \(H\).
	
	\Needspace{5\baselineskip}
	\begin{lemma}\label{lem:coordinate-count}
		For every \(q\in H\),
		\[
		\#\bigl(\overline{L_q}\cap\Sing(\J_d)\bigr)
		\geq\left\lfloor\frac{d+1}{2}\right\rfloor.
		\]
		More precisely, for every \(q\in H\) there is a block
		\(I(q)=\{a+r:0\leq r<\lfloor(d+1)/2\rfloor\}\) of consecutive
		indices modulo \(N\) such that
		\[
		L_q\cap D_{\xi^j}(U_d)\ne\varnothing,\qquad j\in I(q).
		\]
		The corresponding singularities \(p_j\) belong to
		\(\overline{L_q}\). If \(q=[0:0:1]\) or \(q=[1:0:0]\),
		the intersection holds for every index \(j\), and the closure
		contains all \(N\) singular points.
	\end{lemma}
	
	\begin{proof}
		Suppose first that \(q=[u:0:1]\) with \(u\ne0\). In the proof
		of Lemma~\ref{lem:coordinate-line}, the nearest-root condition
		was used only to ensure that the argument of \(\zeta u\) lies in
		\((-\pi/(2d),\pi/(2d))\). For every \(N\)-th root of
		unity satisfying
		\begin{equation}\label{eq:admissible-rotations}
			|\arg(\zeta u)|<\frac{\pi}{2d},
		\end{equation}
		where \(\arg\) is chosen in that interval, the same
		short-time argument shows that the leaf through
		\(D_\zeta(q)\) contains points of \(U_d\). Hence
		the leaf \(L_q\) meets \(D_\zeta^{-1}(U_d)\).
		Proposition~\ref{prop:sector} gives
		\(p_*\in\overline{D_\zeta(L_q)}\), and hence
		\[
		D_\zeta^{-1}(p_*)\in\overline{L_q}.
		\]
		Thus each root satisfying \eqref{eq:admissible-rotations}
		gives a singular point in the closure. More precisely, if
		\(\zeta=\xi^j\), with \(\xi=e^{2\pi\mathrm i/N}\) fixed above, the
		singularity is \([\xi^{-j}:\xi^{dj}:1]\). Distinct roots yield
		distinct singularities by their first affine coordinate.
		
		When \(u\ne0\), the argument of \(u\) determines the position
		of the admissible arc relative to the \(N\)-th roots of unity.
		Its exact number of roots may vary, but the minimum depends
		only on the length of the arc. The admissible rotations
		lie in an open arc of length
		\(\pi/d\) on the unit circle. Since
		\(\xi=e^{2\pi\mathrm i/N}\), the admissible powers \(\xi^j\)
		form a consecutive block of residues modulo \(N\). Moreover,
		\(D_{\xi^j}^{-1}(p_*)=p_{-j}\), so the corresponding singular
		points also have consecutive indices, in reverse order. As consecutive \(N\)-th
		roots of unity are separated by an angle \(2\pi/N\),
		rescaling the angular coordinate reduces the count to the
		integers in an open interval of length
		\[
		\frac{\pi/d}{2\pi/N}
		=\frac{N}{2d}
		=\frac{d+1}{2}+\frac{1}{2d}.
		\]
		Write this length as \(k+\varepsilon\), where \(k\in\mathbb Z_{\geq0}\)
		and \(0<\varepsilon<1\). Irrespective of the position
		of its endpoints, every interval
		\((\alpha,\alpha+k+\varepsilon)\) contains the \(k\) integers
		\(\lfloor\alpha\rfloor+1,\ldots,\lfloor\alpha\rfloor+k\).
		Thus every arc of the stated length contains at least
		\[
		k=\left\lfloor\frac{N}{2d}\right\rfloor
		=\left\lfloor\frac{d+1}{2}\right\rfloor
		\]
		roots. This argument also takes account of the strict
		inequality at the ends of the arc; no genericity assumption
		on the argument of \(u\) is needed.
		
		If \(u=0\), the positive real curve in the proof of
		Lemma~\ref{lem:coordinate-line} shows that \(L_q\) meets
		\(U_d\). Each diagonal symmetry fixes the regular point
		\(q=[0:0:1]\), so it preserves \(L_q\). Hence \(L_q\)
		meets every \(D_{\xi^j}(U_d)\), giving all \(p_j\) in its
		closure. The cyclic symmetry takes \([1:0:0]\) to
		\([0:0:1]\) and preserves \(U_d\), since the sector condition
		is symmetric in the three coordinates. Thus the same argument
		applies to the other coordinate point of \(H\).
	\end{proof}
	
	The cyclic symmetry transfers the same accumulation statement
	to the other coordinate lines, but the three resulting sets of
	singularities need not be disjoint. Their pairwise intersections
	are controlled by the following elementary congruence argument.
	
	\Needspace{5\baselineskip}
	\begin{lemma}\label{lem:progressions}
		Let \(d\geq2\), put \(N=d^2+d+1\), and set
		\(k=\lfloor(d+1)/2\rfloor\). For integers \(a_0,a_1,a_2\),
		consider the subsets of \(\mathbb Z/N\mathbb Z\)
		\[
		I_i=\bigl\{d^{2i}(a_i+r):0\leq r<k\bigr\},
		\qquad i=0,1,2.
		\]
		Then \(\#(I_i\cap I_j)\leq1\) whenever \(i\ne j\), and
		$
		\#(I_0\cup I_1\cup I_2)\geq\max\{k,3k-3\}.
		$
	\end{lemma}
	
	\begin{proof}
		We have \(d^3-1=(d-1)N\), so \(d^3\equiv1\pmod N\),
		and multiplication by \(d\) is invertible modulo \(N\).
		For any two distinct indices we may take them in cyclic order
		as \(i\) and \(i+1\), with indices read modulo \(3\).
		The steps in \(I_0,I_1,I_2\) are \(1,d^2,d^4\), and
		\(d^4\equiv d\pmod N\). For the pairs \((0,1)\), \((1,2)\)
		and \((2,0)\), multiplication by the inverse of the first
		step therefore reduces the two progressions to steps
		\(1,d^2\), since \(d^{-1}\equiv d^2\pmod N\).
		The initial points of the progressions may be arbitrary
		residues; they disappear when differences are taken.
		Suppose there are two common elements, represented by pairs
		\((r,s)\) and \((r',s')\), with \(0\leq r,r',s,s'<k\).
		Taking their differences gives
		\[
		r-r'\equiv d^2(s-s')\pmod N,
		\qquad d(r-r')\equiv s-s'\pmod N.
		\]
		But
		\[
		|d(r-r')-(s-s')|
		\leq(d+1)(k-1)<\frac{d(d+1)}2<N.
		\]
		Hence the last congruence is an equality of integers.
		As \(|s-s'|\leq k-1<d\), it follows that
		\(r=r'\) and \(s=s'\). Thus the two common elements
		coincide. The map \(r\mapsto d^{2i}(a_i+r)\) is injective on
		\(\{0,\ldots,k-1\}\): multiplication by \(d\) is a
		unit modulo \(N\), and \(k<N\). Thus each \(I_i\)
		has exactly \(k\) elements. The three pairwise intersections
		have cardinality at most one, while the triple intersection
		can only increase the cardinality computed by
		inclusion--exclusion. Consequently,
		\[
		\#(I_0\cup I_1\cup I_2)
		\geq\sum_{i=0}^2\#I_i-\sum_{0\leq i<j\leq2}
		\#(I_i\cap I_j)\geq3k-3.
		\]
		The union also contains each \(I_i\), which has cardinality
		\(k\), proving the stated maximum.
	\end{proof}
	
	For \(d=3\), the three sets \(\{0,1\}\), \(\{4,0\}\)
	and \(\{1,4\}\) modulo \(13\) have steps \(1,9,3\),
	respectively. Their three pairwise intersections are distinct
	and their triple intersection is empty. Thus the pairwise
	overlap information alone cannot improve \(3k-3\).
	
	\begin{proof}[Proof of Theorem~\ref{thm:main}]
		Apply Theorem~\ref{thm:propagation-count} to
		\(X=\PP^2\), \(\mathcal F=\J_d\), and the singular
		points \(p_j\) with attracting domains
		\(U_{p_j}=D_{\xi^j}(U_d)\).
		Proposition~\ref{prop:sector} and pushforward supply the
		required complete positive flows. The domain \(U_d\) is
		invariant under cyclic permutation, so cyclic images give
		the same chosen domains after reindexing.
		
		The three curves are \(H_i=A^i(H)\), \(i=0,1,2\).
		Their complements are affine planes, and they avoid the
		singular set by \eqref{eq:singular-set}. Set
		\(k=\lfloor(d+1)/2\rfloor\). For \(q\in H_i\), put
		\(q'=A^{-i}(q)\in H\). Lemma~\ref{lem:coordinate-count}
		gives a consecutive block \(\{a+r:0\leq r<k\}\) with
		\[
		L_{q'}\cap D_{\xi^{a+r}}(U_d)\ne\varnothing.
		\]
		Applying \(A^i\) shows that
		\[
		B_i(q)=\{A^i(p_{a+r}):0\leq r<k\}\subset A_{U_d}(q).
		\]
		This includes coordinate points, for which every index is
		admissible. By \eqref{eq:cyclic-indices}, the index set of
		\(B_i(q)\) is
		\[
		\{d^{2i}(a+r):0\leq r<k\}.
		\]
		The entry calculation for each of the finitely many singular points
		in \(B_i(q)\) has a strict angular margin. Fixing small
		entry times makes these paths valid on a neighbourhood of
		\(q\). As explained after the proof of
		Theorem~\ref{thm:propagation-count}, the resulting local
		paths give its holonomy alternative on that neighbourhood,
		with endpoints inside the attracting domains. Compactness
		of \(H_i\) gives a finite cover with these fixed sets of singular points.
		This also handles the coordinate points, where every index
		is available. Lemma~\ref{lem:progressions} bounds the overlap
		of any two blocks of consecutive indices from different lines by one,
		independently of their initial indices. Thus
		Theorem~\ref{thm:propagation-count}, with \(r_i=k\) and
		\(b_{ij}=1\), gives
		\[
		\#\bigl(\overline L\cap\Sing(\J_d)\bigr)
		\geq\max\{k,3k-3\}=b_d.
		\]
		For \(d=2\), this maximum is one; for \(d\geq3\), it
		is \(3\lfloor(d+1)/2\rfloor-3\). Lemma~\ref{lem:zorn}
		excludes exceptional minimal sets.
	\end{proof}
	
	The leaf through a coordinate point belongs to an open family
	of leaves with the same singular accumulation property. This
	and the numerical bound persist by the general theorem of
	Section~\ref{sec:persistence}.
	
	\Needspace{5\baselineskip}
	\begin{corollary}\label{cor:robustness}
		For each \(d\geq2\), there is a neighbourhood
		\(\mathcal V_d\) of \(\J_d\) in
		\(\Fol(\PP^2,d)=\Fol(\PP^2,\mathcal O_{\PP^2}(d+2))\)
		such that every regular leaf \(L\) of every
		\(\mathcal F\in\mathcal V_d\) satisfies
		\[
		\#\bigl(\overline L\cap\Sing(\mathcal F)\bigr)\geq b_d.
		\]
		For every such \(\mathcal F\), there is also a nonempty
		open saturated subset \(\mathcal A_{\mathcal F}\) of the
		regular locus such that
		\(\Sing(\mathcal F)\subset\overline L\) for every leaf
		\(L\subset\mathcal A_{\mathcal F}\).
		In particular, these foliations have no exceptional minimal set.
	\end{corollary}
	
	\begin{proof}
		All singularities of \(\J_d\) are hyperbolic by
		Remark~\ref{rem:eigenvalues}. Apply
		Theorem~\ref{thm:persistence} with
		\(P=\Sing(\J_d)\), \(D=H\), and the uniform bound
		\(b_d\) already proved in Theorem~\ref{thm:main}.
		There are exactly \(N=d^2+d+1\) continued singularities:
		the defining section is bounded away from zero outside
		small neighbourhoods of the original zeros, and the implicit
		function theorem gives exactly one zero in each of them
		after shrinking the parameter neighbourhood.
		
		For \(\mathcal F=\J_d\), let
		\[
		\mathcal A_j=\{q\in\PP^2\setminus\Sing(\J_d):
		L_q\cap D_{\xi^j}(U_d)\ne\varnothing\}.
		\]
		It is the saturation of the open regular part of
		\(D_{\xi^j}(U_d)\), hence is open and saturated.
		Lemma~\ref{lem:coordinate-count} places
		\(q_0=[0:0:1]\) in every \(\mathcal A_j\).
		Thus \(\bigcap_{j=0}^{N-1}\mathcal A_j\) is nonempty,
		open and saturated, and every leaf in it accumulates at all
		singularities. The final assertion of
		Theorem~\ref{thm:persistence}, applied to \(q_0\), gives
		the corresponding open set for every nearby foliation.
		Lemma~\ref{lem:zorn} gives the last conclusion.
	\end{proof}
	
	\begin{remark}
		Theorem~\ref{thm:main} gives a uniform lower bound for the
		number of singularities in the closure of any regular leaf.
		It does not establish the optimality of that bound. For the
		leaf through
		\([0:0:1]\), the full singular set belongs to its closure
		by Lemma~\ref{lem:coordinate-count}, as already observed in
		\cite[Proposition~1]{CF}. A diagonal symmetry generally moves
		one leaf to another, so the corresponding conclusion cannot
		be deduced for an arbitrary leaf from symmetry alone.
	\end{remark}
	
	\section{Finite pullbacks}\label{sec:finite}
	
	The preceding argument constructs singular limits on the
	projective plane. A finite morphism transfers them to another
	surface through the behaviour of leaf closures. We first prove
	that behaviour for an arbitrary foliation.
	
	Let \(f:Y\to X\) be a finite surjective holomorphic map
	between compact complex surfaces, and let \(\mathcal F\) be
	a holomorphic foliation on \(X\) with finite singular set
	\(S\). The \emph{saturated pullback}
	\(\mathcal G=f^*\mathcal F\) is obtained by pulling back
	local defining holomorphic one-forms and dividing out their
	common divisorial factors. The resulting local forms have
	only isolated zeros and define a foliation on \(Y\). Its
	singular set is therefore finite.
	
	At a regular point of \(\mathcal F\), let \(h\) be a local
	holomorphic first integral. Then \(h\circ f\) is a first
	integral of \(\mathcal G\), although its differential may
	vanish along a ramification divisor. Saturation removes this
	divisorial vanishing from the defining one-form; it preserves
	the fact that \(h\circ f\) is constant on the leaves.
	
	\Needspace{5\baselineskip}
	\begin{theorem}\label{thm:finite-pullback}
		Suppose \(f\) is unramified over \(S\). For every regular
		leaf \(L\) of \(\mathcal G=f^*\mathcal F\), there is a
		regular leaf \(M\) of \(\mathcal F\) such that
		\begin{equation}\label{eq:cofinite-leaf-image}
			M\setminus f(\Sing(\mathcal G))\subset f(L)\subset M,
			\qquad f(\overline L)=\overline M.
		\end{equation}
		Let \(e_L\) be the degree over \(M\) of the normalisation
		of the irreducible component of \(Y\times_XM\) determined
		by the graph of the leaf. Then \(1\leq e_L\leq\deg f\), and
		\begin{equation}\label{eq:finite-singular-count}
			\#\bigl(\overline L\cap\Sing(\mathcal G)\bigr)
			\geq\#\bigl(\overline L\cap f^{-1}(S)\bigr)
			\geq e_L\,\#\bigl(\overline M\cap S\bigr).
		\end{equation}
		In particular, every uniform lower bound for singular points
		in regular leaf closures of \(\mathcal F\) also holds for
		\(\mathcal G\).
	\end{theorem}
	
	\begin{proof}
		Over each point of \(S\), the map \(f\) is a local
		biholomorphism. Hence
		$
		f^{-1}(S)\subset\Sing(\mathcal G).
		$
		The image of \(L\) therefore lies in the regular locus of
		\(\mathcal F\). Local first integrals are constant on its
		image, so connectedness places \(f(L)\) in one regular leaf
		\(M\). A finite map cannot send a curve to a point; thus
		\(f|_L\) is nonconstant.
		
		Give \(M\) its intrinsic Riemann surface structure, and write
		\(\iota:M\to X\) for the leaf inclusion. Form the analytic
		fibre product
		\[
		C=Y\mathbin{\times}_X M.
		\]
		The projection \(C\to M\) is finite, since finiteness is
		preserved by base change. Choose a point \(q\in L\) at which
		\(f|_L\) is noncritical and whose image avoids
		\(f(\Sing(\mathcal G))\). Such a point exists because the
		map is nonconstant and the latter set is finite. The graph
		of the local map from \(L\) to \(M\) gives a curve branch
		of \(C\) at \((q,\iota^{-1}(f(q)))\). Let \(C_0\) be an
		irreducible component containing this branch, endowed with
		its reduced structure.
		
		The map from \(L\) to the intrinsic Riemann surface \(M\)
		is holomorphic in local plaque coordinates. Its graph defines
		a holomorphic map \(L\to C\). The inverse image of the
		closed analytic subset \(C_0\) contains a neighbourhood of
		\(q\), so the identity theorem on the connected leaf puts
		the whole graph in \(C_0\). The projection of \(C_0\)
		to \(M\) is nonconstant and finite. Its image is a closed analytic subset of the
		connected Riemann surface \(M\), and has dimension one;
		hence it is all of \(M\). Normalise \(C_0\). We obtain a
		connected Riemann surface \(\widetilde C\), a finite
		surjective map \(\pi:\widetilde C\to M\) of degree \(e_L\),
		and a map
		\(j:\widetilde C\to Y\) satisfying
		\[
		f\circ j=\iota\circ\pi.
		\]
		
		The set \(T=j^{-1}(\Sing(\mathcal G))\) is finite. Indeed,
		for each singular point \(s\) of \(\mathcal G\), the fibre
		product has at most the point \((s,\iota^{-1}(f(s)))\)
		above \(s\), because the leaf inclusion is injective.
		Normalisation has finite fibres. Removing \(T\) leaves
		\(\widetilde C\setminus T\) connected.
		
		We claim that \(j(\widetilde C\setminus T)\subset L\).
		To see this locally, choose a first integral \(h\) for
		\(\mathcal F\) near \(f(j(z))\). On a small neighbourhood
		of \(z\) in \(\widetilde C\), the function
		\(h\circ f\circ j\) is constant. At the regular point
		\(j(z)\), choose flow-box coordinates for \(\mathcal G\)
		with transverse coordinate \(t\). Since \(h\circ f\) is
		a first integral, it is of the form \(\phi(t)\). The
		function \(\phi\) is nonconstant: a finite map is locally
		dominant, so the pullback of the submersion \(h\) is not
		constant on an open set. The level set of \(\phi\) consists
		locally of isolated transverse values. Thus the local image
		of \(j\) lies in a single plaque of \(\mathcal G\).
		Connectedness of \(\widetilde C\setminus T\), and the branch
		through \(q\in L\), prove the claim.
		
		Now let \(x\in M\setminus f(\Sing(\mathcal G))\). Choose
		\(z\in\pi^{-1}(x)\). Then \(j(z)\) is regular for
		\(\mathcal G\), so \(z\notin T\) and \(j(z)\in L\).
		This proves the first inclusion in
		\eqref{eq:cofinite-leaf-image}. Removing finitely many points
		from a Riemann surface leaves a dense subset, so
		\(\overline{f(L)}=\overline M\). Since \(Y\) is compact,
		\(f(\overline L)=\overline{f(L)}\), which gives the asserted
		equality of closures.
		
		Fix \(s\in\overline M\cap S\). Choose a small connected
		neighbourhood \(V\) of \(s\), disjoint from the branch
		locus, whose inverse image is a disjoint union of open sets
		mapped biholomorphically to \(V\). Take a sequence
		\(x_n\in M\cap V\), tending to \(s\) and avoiding
		\(f(\Sing(\mathcal G))\). Such a sequence exists because
		this excluded set is finite and \(s\) is a limit of the
		regular leaf.
		
		Above \(M\cap V\), the fibre product is unramified, so its
		normalisation is unchanged and \(\pi^{-1}(x_n)\) has
		\(e_L\) distinct points. Their images under \(j\) are
		distinct points of \(L\), by the claim above, and occupy
		distinct local sheets of \(f\). There are only finitely
		many choices of \(e_L\) sheets. Passing to a subsequence
		keeps that choice fixed. The inverse branches of \(f\)
		then show that these \(e_L\) points tend to \(e_L\)
		distinct points of \(\overline L\cap f^{-1}(s)\).
		All are singular for \(\mathcal G\). Fibres over different
		points of \(S\) are disjoint, which proves
		\eqref{eq:finite-singular-count}.
	\end{proof}
	
	The finite exceptional set in \eqref{eq:cofinite-leaf-image}
	accounts for singularities created above tangencies to the
	branch locus. For example, take the local map
	\(f(u,v)=(u,v^2)\) and the regular foliation with first
	integral \(y-x^2\). Its pullback has first integral
	\(v^2-u^2\) and is singular at the origin. The two leaves
	\(v=\pm u\), with the origin removed, each map to the
	parabola \(y=x^2\) with its origin removed. The equality
	of closures in \eqref{eq:cofinite-leaf-image} retains that
	point.
	
	If a base leaf \(M\) is not contained in the branch locus,
	the inverse image of its complement is a finite unramified
	cover of the corresponding open Riemann surface. When the
	monodromy of this cover is transitive, its normalisation is
	connected and \(e_L=\deg f\). Without this condition the
	degrees of the several components can be smaller; the degree
	of \(f\) alone need not multiply the bound for each lifted
	leaf.
	
	\Needspace{5\baselineskip}
	\begin{corollary}\label{cor:finite-jouanolou}
		Let \(Y\) be a smooth projective complex surface and
		\(f:Y\to\PP^2\) a finite surjective morphism unramified
		over \(\Sing(\J_d)\). Then the saturated pullback
		\(\mathcal G=f^*\J_d\) has no invariant algebraic curve,
		and every regular leaf \(L\) satisfies
		\[
		\#\bigl(\overline L\cap\Sing(\mathcal G)\bigr)\geq b_d.
		\]
	\end{corollary}
	
	\begin{proof}
		The inequality follows from Theorems~\ref{thm:finite-pullback}
		and \ref{thm:main}. If \(C\subset Y\) were an invariant
		algebraic curve, its image \(f(C)\) would be an algebraic
		curve, since \(f\) is finite. At a general point of \(C\),
		the restriction of \(f\) to the normalisation of \(C\) has
		nonzero derivative and sends the tangent line of
		\(\mathcal G\) to that of \(\J_d\). Thus \(f(C)\) would
		be invariant by \(\J_d\), contrary to Jouanolou's theorem.
	\end{proof}
	
	\Needspace{5\baselineskip}
	\begin{proposition}\label{prop:pullback-bundles}
		Let \(f:Y\to\PP^2\) be a finite surjective morphism from a
		smooth projective surface, unramified over \(\Sing(\J_d)\),
		and put \(A=f^*\mathcal O_{\PP^2}(1)\). For
		\(\mathcal G=f^*\J_d\),
		\begin{equation*}
			N_{\mathcal G}\simeq A^{\otimes(d+2)},\qquad
			K_{\mathcal G}\simeq K_Y\otimes A^{\otimes(d+2)},
		\end{equation*}
		and
		\begin{equation}\label{eq:milnor-surface}
			\sum_{p\in\Sing(\mathcal G)}\mu_p(\mathcal G)
			=c_2(Y)+(d+2)K_Y\cdot A+(d+2)^2A^2.
		\end{equation}
		Here \(K_{\mathcal G}=\det(T_{\mathcal G})^*\), and
		\(\mu_p\) is the local colength of the coefficients of a
		holomorphic generator.
	\end{proposition}
	
	\begin{proof}
		The normal bundle of a degree \(d\) foliation on \(\PP^2\)
		is \(\mathcal O_{\PP^2}(d+2)\). We check that saturation
		does not remove a divisor from its pullback. Away from
		ramification the map is locally biholomorphic, and the inverse
		image of the singular set is finite. At a general point of a
		ramification component, choose coordinates in which
		\(f(u,v)=(u,v^e)\). If the local defining form on the base is
		\(a(x,y)\,dx+b(x,y)\,dy\), its pullback is
		\[
		a(u,v^e)\,du+e v^{e-1}b(u,v^e)\,dv.
		\]
		A common factor \(v\) would imply \(a(u,0)\equiv0\),
		so the image branch component would be invariant by \(\J_d\).
		There is no invariant algebraic curve. Thus no common
		divisorial factor occurs, proving the first isomorphism.
		The identity \(K_{\mathcal G}=K_Y\otimes N_{\mathcal G}\)
		gives the second.
		
		The pulled-back form is a section of
		\(\Omega_Y^1\otimes A^{\otimes(d+2)}\) with isolated zeros.
		The length of its zero scheme is its second Chern number.
		For a rank-two bundle \(E\) and a line bundle \(L\),
		\(c_2(E\otimes L)=c_2(E)+c_1(E)c_1(L)+c_1(L)^2\).
		Apply this with \(E=\Omega_Y^1\), using
		\(c_1(\Omega_Y^1)=K_Y\) and \(c_2(\Omega_Y^1)=c_2(Y)\).
		The local lengths equal the Milnor numbers, since contracting
		a generator with a local nonvanishing two-form interchanges
		its two coefficients up to sign.
	\end{proof}
	
	\begin{proof}[Proof of Corollary~\ref{cor:any-surface}]
		Choose a very ample line bundle \(A\) on \(Y\). Three general
		sections have no common zero and define a morphism
		\(h:Y\to\PP^2\) with \(h^*\mathcal O_{\PP^2}(1)=A\).
		A curve in a fibre would have degree zero with respect to
		\(A\), contrary to ampleness. Thus the fibres are finite;
		properness makes \(h\) finite, and its two-dimensional image
		is all of \(\PP^2\). This is also the construction by finite
		linear projection; see \cite[Tag~0B1N]{Stacks}.
		
		Choose \(d\) such that \(b_d\geq r\), and let \(B\) be the
		branch locus of \(h\). For each \(s\in\Sing(\J_d)\), the
		condition \(s\notin a(B)\) defines a nonempty Zariski open
		subset of \(\operatorname{PGL}_3(\C)\). Their finite
		intersection is nonempty, since this group is irreducible.
		Choose \(a\) in it and set \(f=a\circ h\).
		Then \(f\) is unramified over \(\Sing(\J_d)\).
		
		By Corollary~\ref{cor:finite-jouanolou},
		\(\mathcal G=f^*\J_d\) has no invariant algebraic curve.
		Theorem~\ref{thm:finite-pullback} gives at least \(b_d\)
		points of \(P=f^{-1}(\Sing(\J_d))\) in every regular leaf
		closure. Every point of \(P\) is hyperbolic, because \(f\)
		is locally biholomorphic there. Choose a smooth ample divisor
		\(D\) avoiding the finite set \(\Sing(\mathcal G)\).
		Its complement is Stein. Proposition~\ref{prop:pullback-bundles}
		gives \(N_{\mathcal G}=N_r:=A^{\otimes(d+2)}\), and
		Theorem~\ref{thm:persistence} supplies the desired nonempty
		analytic open subset of \(\Fol(Y,N_r)\) containing
		\(\mathcal G\). Lemma~\ref{lem:zorn} excludes compact
		saturated subsets of the regular locus for every foliation
		in that open set.
	\end{proof}
	
	The line bundle \(N_r\) is fixed within each such open set,
	but varies as the requested lower bound increases. For a
	fixed normal bundle \(N\), the total Milnor number is the
	fixed integer \(c_2(Y)+K_Y\cdot N+N^2\), which bounds the
	number of distinct singular points. The construction proves
	the absence of invariant algebraic curves for the displayed
	pullback foliation; it does not require this property for
	every member of the analytic open set.
	
	\Needspace{5\baselineskip}
	\begin{corollary}\label{cor:finite-perturbations}
		Fix \(d\geq2\) and a finite morphism \(f:Y\to\PP^2\)
		unramified over \(\Sing(\J_d)\). Put
		\(N=f^*\mathcal O_{\PP^2}(d+2)\).
		There is a neighbourhood of \(f^*\J_d\) in \(\Fol(Y,N)\)
		such that every regular leaf of every foliation in it has at
		least \(b_d\) singular points in its closure.
	\end{corollary}
	
	\begin{proof}
		Use \(P=f^{-1}(\Sing(\J_d))\) and a smooth ample divisor
		avoiding \(\Sing(f^*\J_d)\). Theorem~\ref{thm:finite-pullback}
		gives the bound using points of \(P\), all of which are
		hyperbolic. Apply Theorem~\ref{thm:persistence} and
		Proposition~\ref{prop:pullback-bundles}.
	\end{proof}
	
	\section{Double planes and K3 surfaces}
	
	Double covers of the plane give examples for which the normal
	bundle, the canonical bundle and all singularities can be
	described explicitly. We use the standard construction of a
	K3 surface as a double plane branched over a smooth sextic;
	see \cite[Chapter~1]{Huybrechts}.
	
	\Needspace{5\baselineskip}
	\begin{lemma}\label{lem:sextic-tangencies}
		For fixed \(d\geq2\), a general plane sextic \(B\) is smooth,
		avoids \(\Sing(\J_d)\), and has only simple tangencies with
		\(\J_d\). Their number is \(6(d+5)\).
	\end{lemma}
	
	\begin{proof}
		Smoothness and avoidance of a finite set hold on nonempty
		Zariski open subsets of the space of sextics. At a regular
		point \(x\), choose a local nonvanishing generator \(V\)
		and a local equation \(F\) of a sextic. A tangency of order
		at least two requires
		\(F(x)=VF(x)=V^2F(x)=0\). Sections of
		\(\mathcal O_{\PP^2}(6)\) generate all two-jets at \(x\).
		In flow-box coordinates for \(V\), these three equations
		specify the constant, first and second coefficients in the
		leaf direction and are independent linear conditions.
		Their incidence over the two-dimensional regular locus has
		dimension at most one less than the space of sextics. Its
		image therefore lies in a proper algebraic subset after
		taking closure. A general smooth sextic avoiding the singular
		set has no such point. Its tangencies are simple.
		
		Restriction of a defining form of \(\J_d\) gives a section
		of \(K_B\otimes N_{\J_d}|_B\), whose zeros are precisely
		the tangencies. This section is nonzero because \(B\) is
		not invariant. A smooth plane sextic has genus ten, so
		\(\deg K_B=18\); also
		\(\deg N_{\J_d}|_B=6(d+2)\). All zeros are simple, and
		their number is \(18+6(d+2)=6(d+5)\).
	\end{proof}
	
	\Needspace{5\baselineskip}
	\begin{theorem}\label{thm:k3}
		Choose a sextic \(B\) as in Lemma~\ref{lem:sextic-tangencies},
		and let \(f:Y\to\PP^2\) be the double cover branched over
		\(B\). Then \(Y\) is a K3 surface. Put
		\(A=f^*\mathcal O_{\PP^2}(1)\), so \(A^2=2\).
		The foliation \(\mathcal G=f^*\J_d\) has no invariant
		algebraic curve and satisfies
		$
		N_{\mathcal G}\simeq K_{\mathcal G}\simeq A^{\otimes(d+2)}.
		$
		Every regular leaf has at least \(b_d\) singular points in
		its closure. The singular set consists of
		\(2(d^2+d+1)\) hyperbolic points above \(\Sing(\J_d)\)
		and \(6(d+5)\) nondegenerate singularities with a holomorphic
		Morse first integral above the tangencies with \(B\).
		All singularities are reduced and have Milnor number one;
		in particular,
		\[
		\#\Sing(\mathcal G)=24+2(d+2)^2.
		\]
		The canonical bundle \(K_{\mathcal G}\) is ample and
		\(\mathcal G\) is of general type. The lower bound \(b_d\)
		holds throughout a neighbourhood of \(\mathcal G\) in
		\(\Fol(Y,A^{\otimes(d+2)})\).
	\end{theorem}
	
	\begin{proof}
		The covering construction gives
		$ f_*\mathcal O_Y=\mathcal O_{\PP^2}\oplus\mathcal O_{\PP^2}(-3),
		\ K_Y=f^*(K_{\PP^2}\otimes\mathcal O_{\PP^2}(3))
		\simeq\mathcal O_Y.
		$
		Thus \(H^1(Y,\mathcal O_Y)=0\).
		Smoothness follows from smoothness of the branch sextic,
		so \(Y\) is a K3 surface. The pullback of an ample line
		bundle under a finite morphism is ample; hence \(A\) is
		ample and \(A^2=\deg f=2\). Proposition~\ref{prop:pullback-bundles}
		gives the asserted bundles, and
		Corollary~\ref{cor:finite-jouanolou} gives the dynamical
		bound and the absence of invariant curves.
		
		Above a singular point of \(\J_d\), the cover is unramified
		and has two sheets, giving the \(2(d^2+d+1)\) hyperbolic
		singularities. At a branch point, write \(f(u,v)=(u,v^2)\)
		with \(B=\{y=0\}\) on the base, and choose a local first
		integral \(h(x,y)\) of the regular foliation. The pullback
		has first integral \(h(u,v^2)\). At a transverse intersection
		with \(B\), \(h_x\ne0\), so the pullback is regular.
		At a simple tangency, \(h_x=0\), \(h_y\ne0\), and
		\(h_{xx}\ne0\). The Hessian of the pulled-back first
		integral is diagonal with nonzero entries \(h_{xx}\) and
		\(2h_y\). The holomorphic Morse lemma gives local
		coordinates in which the first integral is \(v^2-u^2\).
		Its defining form has an isolated zero, its Milnor number is
		one, and a generator has eigenvalue ratio \(-1\).
		
		There are no further singularities: off the ramification and
		the inverse image of \(\Sing(\J_d)\), the map is locally
		biholomorphic. Lemma~\ref{lem:sextic-tangencies} therefore
		gives the asserted list. The hyperbolic points are reduced,
		as are the Morse singularities, whose ratio \(-1\) is not
		a positive rational number. Moreover,
		\[
		2(d^2+d+1)+6(d+5)=24+2(d+2)^2,
		\]
		in agreement with \eqref{eq:milnor-surface}, since
		\(c_2(Y)=24\). This Chern number also follows from
		\(\chi_{\mathrm{top}}(Y)=2\chi_{\mathrm{top}}(\PP^2)
		-\chi_{\mathrm{top}}(B)=6-(2-20)=24\).
		As the foliation is reduced, its Kodaira dimension is the
		Iitaka dimension of \(K_{\mathcal G}\); ampleness makes it
		two, so the foliation is of general type in the usual sense
		\cite{Brunella}. Finally, Corollary~\ref{cor:finite-perturbations}
		gives the neighbourhood asserted in the theorem. It counts
		the continued hyperbolic points above \(\Sing(\J_d)\);
		the Morse singularities are not hyperbolic.
	\end{proof}
	
	\Needspace{5\baselineskip}
	\begin{corollary}
		In Theorem~\ref{thm:k3}, suppose the base leaf \(M\)
		associated with a regular leaf \(L\) meets \(B\)
		transversely at some point. Then \(e_L=2\) and
		\[
		\#\bigl(\overline L\cap f^{-1}(\Sing(\J_d))\bigr)
		\geq2\,\#\bigl(\overline M\cap\Sing(\J_d)\bigr)
		\geq2b_d.
		\]
	\end{corollary}
	
	\begin{proof}
		At a transverse intersection, the restriction of the covering
		to the leaf has local equation \(v^2=t\), where \(t\) is
		a coordinate on \(M\). Its normalisation therefore has a
		ramification point of degree two. A disconnected degree-two
		normalisation would have two components of degree one, which
		cannot have such a point. Hence the normalised fibre product
		is connected and has degree two. Apply
		Theorem~\ref{thm:finite-pullback}.
	\end{proof}

	\subsection*{Acknowledgements.}
	MC acknowledges partial support from the Università degli Studi di Bari Aldo Moro. He is a member of GNSAGA, a national research group of INdAM. 
	We are grateful to Arturo Fernández-Pérez, Fernando Lourenço, Alan Muniz, Miguel Rodríguez Peña, Bruno Scardua, José Seade, and Márcio Soares for many stimulating and illuminating discussions.
	
	\smallskip
	\subsection*{AI disclosure} The mathematical ideas, results and proofs presented in this paper were conceived and developed by the authors. ChatGPT (OpenAI) was used as an auxiliary tool to assist in checking calculations, examining mathematical arguments, and drafting and revising the exposition. The authors take full responsibility for the mathematical correctness and content of the paper.

\end{document}